\documentclass[11pt]{article}

\usepackage{amsmath,amsthm}

\usepackage{amssymb,latexsym}

\usepackage{enumerate}

\newtheorem{tw}{Theorem}[subsection]
\newtheorem{lm}[tw]{Lemma}

\newtheorem{stw}[tw]{Proposition}
\newenvironment{dow}{\it Proof.\rm}{\hfill $\Box$}

\theoremstyle{definition}
\newtheorem*{df}{Definition}
\newtheorem{uw}[tw]{Remark}
\newtheorem{prz}[tw]{Example}

\newcommand{\BL}{{\mathbb L}}

\newcommand{\BN}{{\mathbb N}}
\newcommand{\BR}{{\mathbb R}}

\newcommand{\FF}{{\mathcal{F}}}

\newcommand{\HH}{{\mathcal{H}}}

\newcommand{\MM}{{\mathcal{M}}}

\newcommand{\intt}{{\int_{t}^{T}}}

\newcommand{\BRD}{{\mathbb{R}^{d}}}

\newcommand{\esssup}{\mathop{\mathrm{ess\,sup}}}

\newcommand{\nsubsection}{\setcounter{equation}{0}\subsection}

\begin{document}

\title {Reflected BSDEs with monotone generator}
\author {Tomasz Klimsiak}
\date{}
\maketitle
\begin{abstract}
We give necessary and sufficient condition for existence and
uniqueness of $\mathbb{L}^{p}$-solutions of reflected BSDEs with
continuous barrier, generator monotone with respect to $y$ and
Lipschitz continuous with respect to $z$, and with data in
$\mathbb{L}^{p}$, $p\ge 1$. We also prove that the solutions may
be approximated by the penalization method.
\end{abstract}


\footnotetext{{\em Key words or phrases:} Reflected backward
stochastic differential equation, monotone generator,
$\mathbb{L}^{p}$-solutions.}

\footnotetext{Research supported by the Polish Minister of Science
and Higher Education under Grant N N201 372 436.}

\nsubsection{Introduction}

Let $B$ be a standard $d$-dimensional Brownian motion defined on
some probability space $(\Omega,\FF,P)$ and let $\{\FF_t\}$ denote
the augmentation of the natural filtration generated by $B$. In
the present paper we study the problem of existence, uniqueness
and approximation of $\BL^p$-solutions of reflected backward
stochastic differential equations (RBSDEs for short) with monotone
generator of the form
\begin{equation}
\label{eq1.1} \left\{
\begin{array}{l} Y_{t}=\xi+\intt
f(s,Y_{s},Z_{s})\,ds
-\intt dK_{s} -\intt Z_{s}\,dB_{s},\quad t\in [0,T],\\
Y_{t}\ge L_{t},\quad t\in [0,T], \\
K\mbox{ is continuous, increasing, }K_{0}=0,\,\int_{0}^{T}
(Y_{t}-L_{t})\,dK_{t}=0.
\end{array}
\right.
\end{equation}
Here  $\xi$ is an $\FF_T$-measurable random variable called the
terminal condition, %the random function
$f:[0,T]\times\Omega\times\BR\times\BR^d\rightarrow\BR$ is the
generator (or coefficient) of the equation and an
$\{\FF_t\}$-adapted continuous proces $L=\{L_t,t\in[0,T]\}$  such
that $L_T\le\xi$ $P$-a.s. is called the obstacle (or barrier). A
solution of (\ref{eq1.1}) is a triple $(Y,Z,K)$ of
$\{\FF_t\}$-progressively measurable processes having some
integrability properties depending on assumptions imposed on the
data  $\xi,f,L$ and satisfying (\ref{eq1.1}) $P$-a.s.

Equations of the form (\ref{eq1.1}) were introduced in El Karoui
et al. \cite{EKPPQ}. At present it is widely recognized that they
provide a useful and efficient tool for studying problems in
different mathematical fields, such as mathematical finance,
stochastic control and game theory, partial differential equations
and others (see, e.g., \cite{CK,EKPPQ,EPQ,H,Kl2}).

In \cite{EKPPQ} existence and uniqueness of square-integrable
solutions of (\ref{eq1.1}) are proved under the assumption  that
$\xi$, $\int^T_0|f(t,0,0)|\,dt$ and $L^*_T=\sup_{t\le T}|L_t|$ are
square-integrable, $f$ satisfies the linear growth condition and
is Lipschitz continuous with respect to both variables $y$ and
$z$. These assumptions are too strong for many interesting
applications. Therefore  many attempts have been made to prove
existence and uniqueness of solutions of RBSDEs under less
restrictive assumptions on the data. Roughly speaking one can
distinguish here two types of results: for RBSDEs with less
regular barriers (see, e.g., \cite{PengXu}) and for equations with
continuous barriers whose generators or terminal conditions
satisfy weaker assumptions than in \cite{EKPPQ}. We are interested
in the second direction of investigation of (\ref{eq1.1}).

In the paper we consider $\mathbb{L}^{p}$-integrable data with
$p\ge 1$ and we assume that the generator is continuous and
monotone in $y$ and Lipschitz continuous with respect to $z$.
Assumptions of that type were considered in \cite{A,HP,LMX,RS} but
it is worth mentioning that the case where the generator is
monotone and at the same time the data are
$\mathbb{L}^{p}$-integrable for some $p\in [1,2)$ was considered
previously only in \cite{A,RS} (to be exact, in \cite{A} the
author considers the case $p\in (1,2)$ but for generalized
RBSDEs). Let us also mention that in the case $p=2$ existence and
uniqueness results are known  for equations with generators
satisfying even weaker regularity conditions. For instance, in
\cite{C} continuous generators satisfying the linear growth
conditions are considered, in \cite{ZZ} it is assumed that the
generator is left-Lipschitz continuous and possibly discontinuous
in $y$, and in \cite{K} equations with generators satisfying the
superlinear growth condition with respect to $y$,  the quadratic
growth condition with respect to $z$ and with data ensuring
boundedness of the first component $Y$ are considered. In all
these  papers except for \cite{RS} the authors consider the
so-called general growth condition which says that
\begin{align}\label{i4}
|f(t,y,0)|\le |f(t,0,0)|+\varphi(|y|),\quad
t\in[0,T],y\in\mathbb{R},
\end{align}
where $\varphi:\mathbb{R}^{+}\rightarrow \mathbb{R}^{+}$ is a
continuous increasing function or continuous function which is
bounded on bounded subsets of $\mathbb{R}$. In \cite{RS} weaker
than (\ref{i4}) condition of the form
\begin{align}\label{i5}
\forall_{r>0}\quad \sup_{|y|\le r}|f(\cdot,y,0)-f(\cdot,0,0)|\in
\mathbb{L}^{1}(0,T).
\end{align}
is assumed. Condition (\ref{i5}) seems to be the best possible
growth condition on $f$ with respect to $y$. It was used earlier
in the paper \cite{BDHPS} devoted to $\BL^p$-solutions of usual
(non-reflected) BSDEs with monotone generators. Similar condition
is widely used in the theory of partial differential equations
(see \cite{Betal.} and the references given there). Let us point
out, however, that in contrast to the case of usual BSDEs with
monotone generators, in general assumption (\ref{i4}) (or
(\ref{i5})) together with $\mathbb{L}^{p}$-integrability of the
data (integrability of $\xi$, $L^*_T$, $\int^T_0|f(t,0,0)|\,dt$ in
our case) do not guarantee existence of
$\mathbb{L}^{p}$-integrable solutions of (\ref{eq1.1}). For
existence some additional assumptions relating the growth of $f$
with that of the barrier is required. In \cite{A,LMX} existence of
solutions is proved under the assumption that
$E|\varphi(\sup_{t\le T}e^{\mu t} L^{+}_{t})|^2<+\infty$, where
$\varphi$ is the function of condition (\ref{i4}) and $\mu$ is the
monotonicity coefficient of $f$. In \cite{RS} it is shown that it
suffices to assume that
\begin{align}\label{i7}
E(\int_{0}^{T}|f(t,\sup_{s\le t}
L^{+}_{t},0)|\,dt)^{p}\,dt<+\infty.
\end{align}

Condition (\ref{i7}) is still not the best possible. In our main
result of the paper we give a necessary and sufficient condition for
existence and uniqueness of $\mathbb{L}^{p}$-integrable solution
of RBSDE (\ref{eq1.1}) under the assumptions that the data are
$\BL^p$-integrable, $f$ is monotone in $y$ and Lipschitz
continuous in $z$ and (\ref{i5}) is satisfied. Moreover, our
condition is not only weaker than (\ref{i7}) but at the same time
much easier to check than (\ref{i7}) in case of very important in
applications Markov type RBSDEs with obstacles of the form
$L=h(\cdot,X)$, where $h:[0,T]\times\BR^d\rightarrow\BR$ is a
measurable function and $X$ is a Hunt process associated with some
Markov semigroup. In the case of Markov RBSDEs which appear for
instance in applications to variational problems for PDEs (see,
e.g., \cite{EKPPQ,Kl2}) our condition can be formulated in terms
of $f,h$ only. We prove the main result for $p\ge1$. Moreover, we
show that for $p\ge1$ a unique solution of RBSDE (\ref{eq1.1}) can
be approximated via penalization. The last result strengthens the
corresponding result in \cite{RS} proved in case $p>1$ for general
generators and in case $p=1$ for generators not depending on $z$.

In the last part of the paper we study (\ref{eq1.1}) in the case
where $\xi$, $L^{+,*}$, $\int^T_0|f(t,0,0)|\,dt$ are
$\BL^p$-integrable for some $p\ge 1$ but our weaker form of
(\ref{i7}) is not satisfied. We have already mentioned, that then
there are no $\BL^p$-integrable solutions of (\ref{eq1.1}). We
show that still there exist solutions of (\ref{eq1.1}) having
weaker regularity properties.

The paper is organized as follows.  Section \ref{sec2} contains
notation and main hypotheses used in the paper. In Section
\ref{sec3} we show basic a priori estimates for solutions of
BSDEs. In Section \ref{sec4} we prove comparison results as well
as some useful results on c\`adl\`ag regularity of monotone
limits of semimartingales and uniform estimates of monotone
sequences. In Section \ref{sec5} we prove our main existence and
uniqueness result for $p>1$, and in Section \ref{sec6} for $p=1$.
Finally, in Section \ref{sec7} we deal with nonintegrable
solutions.

\nsubsection{Notation and hypotheses} \label{sec2}

Let $B=\{B_{t}, t\ge 0\}$ be a standard $d$-dimensional Brownian
motion defined on some complete probability space $(\Omega,\FF,P)$
and let $\{\FF_{t}, t\ge 0\}$ be the augmented filtration
generated by $B$. In the whole paper all notions whose definitions
are related to some filtration are understood with respect to the
filtration $\{\FF_{t}\}$.

Given a stochastic process $X$ on $[0,T]$ with values in
$\mathbb{R}^{n}$ we set $X^{*}_{t}=\sup_{0\le s\le t}|X_{s}|$,
$t\in[0,T]$, where $|\cdot|$ denotes the Euclidean norm on
$\mathbb{R}^{n}$.  By $\mathcal{S}$ we denote the set of all
progressively measurable continuous processes. For $p>0$ we denote
by $\mathcal{S}^{p}$ the set of all processes $X\in\mathcal{S}$
such that
\[
\|X\|_{\mathcal{S}^{p}} =(E\sup_{t\in [0,T]}
|X_{t}|^{p})^{1\wedge1/p}<+\infty.
\]
$M$ is the set of all progressively measurable processes $X$ such
that
\[
P(\int_{0}^{T}|X_{t}|^{2}\,dt<+\infty)=1
\]
and for $p>0$, $M^{p}$ is the set of all processes $X\in M$ such
that
\[
(E(\int_{0}^{T}|X_{t}|^{2}\,dt)^{p/2})^{1\wedge 1/p} <+\infty.
\]
For $p,q>0$, $\mathbb{L}^{p,q}(\FF)$ (resp.
$\mathbb{L}^{p}(\FF_{T})$)  denotes the set of all progressively
measurable processes ($\FF_T$ measurable random variables) $X$
such that
\[
(E(\int_{0}^{T}|X_{t}|^{p}\,dt)^{q/(1\wedge 1/p)})^{1\wedge
1/q}<+\infty \quad \left(\mbox{resp. }(E|X|^{p})^{1/p}<+\infty
\right).
\]
For brevity we denote $\mathbb{L}^{p,p}(\FF)$ by
$\mathbb{L}^{p}(\FF)$. By $\mathbb{L}^{1}(0,T)$ we denote the
space of Lebesgue integrable real valued functions on $[0,T]$.

$\MM_{c}$ is the set of all continuous martingales (resp. local
martingales) and $\MM^{p}_{c}$, $p\ge1$, is the set of all
martingales $M\in\MM_{c}$ such that $E(\langle M
\rangle_{T})^{p/2}<+\infty$. $\mathcal{V}_{c}$ (resp.
$\mathcal{V}^{+}_{c}$) is the set of all continuous progressively
measurable processes of finite variation (resp. increasing
processes) and $\mathcal{V}^{p}_{c}$ (resp.
$\mathcal{V}^{+,p}_{c}$) is the set of all processes
$V\in\mathcal{V}_{c}$ (resp. $V\in\mathcal{V}^{+}_{c}$) such that
$E|V|^{p}_{T}<+\infty$. We put
$\HH^{p}_{c}=\MM_{c}^{p}+\mathcal{V}_{c}^{p}$.

For a given measurable process $Y$ of class (D) we denote
\[
\|Y\|_{1}=\sup\{E|Y_{\tau}|,\tau\in\mathcal{T}\}.
\]

In what follows
$f:[0,T]\times\Omega\times\mathbb{R}\times\BRD\rightarrow
\mathbb{R}$ is a measurable function with respect to
$Prog\times\mathcal{B}(\mathbb{R})\times\mathcal{B}(\BRD)$, where
$Prog$ denotes the $\sigma$-field of progressive subsets of
$[0,T]\times\Omega$.

In the whole paper all equalities and inequalities between random
elements are understood to hold $P$-a.s.

Let $p\ge 1$. In the paper we consider the following hypotheses.

\begin{enumerate}
\item[(H1)] $E|\xi|^{p}+E(\int_{0}^{T}|f(t,0,0)|\,dt)^{p}<\infty$
\item[(H2)]  There exists $\lambda>0$ such that $|f(t,y,z)-f(t,y,z')|
\le \lambda |z-z'|$
for every $t\in [0,T], y\in\mathbb{R}, z,z'\in\BRD$.
\item[(H3)] There exists $\mu\in\mathbb{R}$ such that
$(f(t,y,z)-f(t,y',z))(y-y')\le \mu(y-y')^{2}$
for every $t\in [0,T], y, y'\in\mathbb{R}, z,z'\in\BRD$.
\item [(H4)] For every $(t,z)\in[0,T]\times\BRD$ the mapping
$\mathbb{R}\ni y\rightarrow f(t,y,z)$
is continuous.
\item[(H5)] For every $r>0$ the mapping
$[0,T]\ni t\rightarrow\sup_{|y|\le r}|f(t,y,0)-f(t,0,0)|$ belongs
to  $\mathbb{L}^{1}(0,T)$.
\item[(H6)]$L$ is a continuous, progressively
measurable process such that $L_T\le\xi$.
\item [(H7)]There exists a semimartingale $X$ such that
$X\in \HH^{p}_{c}$ for some $p>1$, $X_{t}\ge L_{t}$, $t\in [0,T]$
and $E(\int_{0}^{T}f^{-}(s,X_{s},0)\,ds)^{p}<\infty$.
\item [(H7*)] There exists a semimartingale $X$ of class (D) such that
$X\in\mathcal{V}^{1}_{c}+\MM^{q}_{c}$ for every $q\in (0,1)$,
$X_{t}\ge L_{t}$, $t\in [0,T]$ and
$E\int_{0}^{T}f^{-}(s,X_{s},0)\,ds<+\infty$.
\item[(A)]There exist $\mu\in\mathbb{R}$ and $\lambda \ge 0$ such
that
\[
\hat{y}f(t,y,z)\le f_{t}+\mu|y|+\lambda |z|
\]
for every $t\in[0,T]$, $y\in\mathbb{R}$, $z\in\mathbb{R}^d$, where
$\hat{y}=\mathbf{1}_{\{y\neq 0\}}\frac{y}{|y|}$ and
$\{f_{t};t\in[0,T\}$ is a nonnegative progressively measurable
process.
\item[(Z)]There exist $\alpha\in(0,1)$, $\gamma\ge 0$ and a nonnegative
process $g\in\mathbb{L}^{1}(\FF)$ such that
\[
|f(t,y,z)-f(t,y,0)|\le \gamma (g_{t}+|y|+|z|)^{\alpha}
\]
for every $t\in[0,T]$, $y\in\mathbb{R}$, $z\in\mathbb{R}^d$.
\end{enumerate}

\nsubsection{A priori estimates}
\label{sec3}

In this section $K$ denotes an arbitrary but fixed process of the
class $\mathcal{V}^{+}_{c}$ such that $K_{0}=0$.

The following version of It\^o's formula will be frequently used
in the paper.

\begin{stw}
\label{prop.ito} Let $p\ge 1$ and let $X$ be a progressively
measurable process of the form
\begin{align*}
X_{t}=X_{0}+\int_{0}^{t} dK_{s}+\int_{0}^{t} Z_{s}\,dB_{s},\quad t\in [0,T],
\end{align*}
where $Z\in M$. Then there is $L\in \mathcal{V}^{+}_{c}$ such that
\begin{align*}
|X_{t}|^{p}-|X_{0}|^{p}&=p\int_{0}^{t}
|X_{s}|^{p-1}\hat{X}_{s}\,dK_{s}
+p\int_{0}^{t}|X_{s}|^{p-1}\hat{X}_{s}\,dB_{s}\nonumber\\
&\quad+c(p)\int^t_0\mathbf{1}_{\{X_{s}\neq
0\}}|X_{s}|^{p-2}|Z_{s}|^{2}\,ds +L_{t}\mathbf{1}_{\{p=1\}}
\end{align*}
with $c(p)=p(p-1)/2$.
\end{stw}
\begin{dow}
\!\!The proof is a matter of slight modification of the proof of
\cite[\!Lemma 2.2]{BDHPS}.
\end{dow}

\begin{df}
We say that a pair $(Y,Z)$ of progressively measurable processes
is a solution of BSDE$(\xi,f+dK)$ iff $Z\in M$, the mapping
$[0,T]\ni t\mapsto f(t,Y_{t},Z_{t})$ belongs to
$\mathbb{L}^{1}(0,T)$, $P$-a.s. and
\begin{equation}
\label{eq3.02} Y_{t}=\xi+\int_{t}^{T}f(s,Y_{s},Z_{s})\,ds
+\int_{t}^{T}dK_{s}-\int_{t}^{T} Z_{s}\,dB_{s},\quad t\in[0,T].
\end{equation}
\end{df}

\begin{lm}\label{lm1}
Let $(Y,Z)$ be a solution of \mbox{\rm BSDE}$(\xi,f+dK)$. Assume
that \mbox{\rm(H3)} is satisfied and there exists a progressively
measurable process $X$ such that $X_{t}\ge Y_{t}$, $t\in [0,T]$
and the mappings $[0,T]\ni t\mapsto X^{+}_t$, $[0,T]\ni t\mapsto
f^{-}(t,X_{t},0)$ belong to $\mathbb{L}^{1}(0,T)$, $P$-a.s.
\begin{enumerate}
\item[\rm(i)]If \mbox{\rm(H2)} is satisfied then for every stopping time
$\tau\le T$  and $a\ge \mu$,
\begin{align*}
\int_{0}^{\tau} e^{at}dK_{t}&\le |e^{a\tau}Y_{\tau}|
+|Y_{0}|+\int_{0}^{\tau}e^{as}Z_{s}\,dB_{s}
+\lambda\int_{0}^{\tau}e^{as}|Z_{s}|\,ds\\&\quad
+\int_{0}^{\tau}e^{as}f^{-}(s,X_{s},0)\,ds
+\int_{0}^{\tau}a^{+}e^{as}X_{s}^{+}\,ds.
\end{align*}
\item[\rm(ii)] If \mbox{\rm(Z)} is satisfied then for every stopping
time $\tau\le T$ and $a\ge\mu$,
\begin{align*}
\int_{0}^{\tau} e^{at}dK_{t}&\le |e^{a\tau}Y_{\tau}|+|Y_{0}|
+\int_{0}^{\tau}e^{as}Z_{s}\,dB_{s}+\gamma\int_{0}^{\tau}e^{as}(g_{s}
+|Y_{s}|+|Z_{s}|)^{\alpha}\,ds\\
&\quad+\int_{0}^{\tau}e^{as}f^{-}(s,X_{s},0)\,ds
+\int_{0}^{\tau}a^{+}e^{as}X_{s}^{+}\,ds.
\end{align*}
\end{enumerate}
\end{lm}
\begin{dow}
Assume that $\mu\le 0$. Then $f^{-}(s,Y_{s},0)\le
f^{-}(s,X_{s},0)$, $s\in [0,T]$ and from (\ref{eq3.02}) and (H2) it follows
that
\begin{align*}
K_{\tau}\le -Y_{\tau}+Y_{0}+\int_{0}^{\tau}Z_{s}\,dB_{s}
+\lambda\int_{0}^{\tau}|Z_{s}|\,ds-\int_{0}^{\tau}
f(s,Y_{s},0)\,ds,
\end{align*}
which implies (i) with $a=0$. Now, let $a\ge\mu$ and let
$\tilde{Y}_{t}=e^{at}Y_{t}$, $\tilde{Z}_{t}=e^{at}Z_{t}$ and
$\tilde{\xi}=e^{aT}\xi$, $\tilde{f}(t,y,z)
=e^{at}f(t,e^{-at}y,e^{-at}z)-ay$,
$d\tilde{K}_{t}=e^{at}\,dK_{t}$. Then $\tilde{f}$ satisfies (H3)
with $\mu=0$ and by It\^o's formula,
\[
\tilde{Y}_{t}=\tilde{\xi}
+\int_{t}^{T}\tilde{f}(s,\tilde{Y}_{s},\tilde{Z}_{s})\,ds
+\int_{t}^{T}d\tilde{K}_{s}-\int_{t}^{T}\tilde{Z}_{s}\,dB_{s},
\quad t\in [0,T],
\]
from which in the same manner as before we obtain (i) for
$a\ge\mu$.

To prove (ii) let us observe that from (\ref{eq3.02}) and (Z) it
follows immediately that
\begin{align*}
K_{\tau}\le -Y_{\tau}+Y_{0}+\int_{0}^{\tau}Z_{s}\,dB_{s}
+\gamma\int_{0}^{\tau}(g_{s}+|Y_{s}|+|Z_{s}|)^{\alpha}\,ds
-\int_{0}^{\tau} f(s,Y_{s},0)\,ds.
\end{align*}
Therefore repeating arguments from the proof of (i) we get (ii).
\end{dow}

\begin{lm}\label{lm2}
Assume \mbox{\rm(A)} and let $(Y,Z)$ be a solution of \mbox{\rm
BSDE}$(\xi,f+dK)$. If $Y\in\mathcal{S}^{p}$ for some $p>0$ and
\[
E(\int_{0}^{T}X^{+}_{s}\,ds)^{p}
+E(\int_{0}^{T}f^{-}(s,X_{s},0)\,ds)^{p}
+E(\int_{0}^{T}|f(s,0,0)|\,ds)^{p}<+\infty
\]
for some progressively measurable process $X$ such that $X_{t}\ge
Y_{t}$, $t\in [0,T]$, then $Z\in M^{p}$ and there exists $C$
depending only on $\lambda,p,T$ such that for every $a\ge
\mu+\lambda^{2}$,
\begin{align*}
&E\bigg((\int_{0}^{T} e^{2as}|Z_{s}|^{2}\,ds)^{p/2}+(\int_{0}^{T}
e^{as}\,dK_{s})^{p}\bigg) \le CE\bigg(\sup_{t\le T}
e^{apt}|Y_{t}|^{p}\\
&\qquad+(\int_{0}^{T}e^{as}|f(s,0,0)|\,ds)^{p}
+(\int_{0}^{T}e^{as}f^{-}(s,X_{s},0)\,ds)^{p} +
(\int_{0}^{T}a^{+}e^{as}X_{s}^{+}\,ds)^{p}\bigg).
\end{align*}
\end{lm}
\begin{dow}
By standard arguments we may assume that $\mu+\lambda^{2}\le 0$
and take $a=0$. For each $k\in\mathbb{N}$ let us consider the
stopping time
\begin{equation}
\label{eq3.01} \tau_k=\inf\{t\in [0,T];\int_0^t|Z_{s}|^{2}\,ds\ge
k\}\wedge T.
\end{equation}
Then as in the proof of Eq. (5) in \cite{BDHPS} we get
\begin{align*}
(\int_{0}^{\tau_{k}}|Z_{s}|^{2}\,ds)^{p/2}\le
c_{p}\bigg(|Y^{*}_{T}|^{p} +(\int_{0}^{T} f_{s}\,ds)^{p}
+|\int_{0}^{\tau_{k}}Y_{s}Z_{s}\,dB_{s}|^{p/2}
+(\int_{0}^{\tau_{k}}|Y_{s}|\,dK_{s})^{p/2}\bigg),
\end{align*}
and hence, repeating arguments following Eq. (5) in \cite{BDHPS}
we show that
\begin{align}\label{eq2.1}
E(\int_{0}^{\tau_{k}}|Z_{s}|^{2}\,ds)^{p/2}\le
c_{p}E\bigg(|Y^{*}_{T}|^{p} +(\int_{0}^{T} f_{s}\,ds)^{p}
+(\int_{0}^{\tau_{k}}|Y_{s}|\,dK_{s})^{p/2}\bigg).
\end{align}
By  Lemma \ref{lm1} and the Burkholder-Davis-Gundy inequality,
\begin{align}\label{eq2.2}
EK^{p}_{\tau_{k}}\le c'(p,\lambda,T)E\{|Y^{*}_{T}|^{p}
+(\int_{0}^{\tau_{k}}|Z_{s}|^{2}\,ds)^{p/2}
+(\int_{0}^{T}f^{-}(s,X_{s},0)\,ds)^{p}\}.
\end{align}
Moreover, applying Young's inequality we conclude from
(\ref{eq2.1}) that for every $\alpha>0$,
\begin{align}\label{eq2.3}
&E(\int_{0}^{\tau_{k}}|Z_{s}|^{2}\,ds)^{p/2} \nonumber\\
&\quad\le c''(p,\alpha)E\{|Y^{*}_{T}|^{p}+(\int_{0}^{T}
f_{s}\,ds)^{p}+(\int_{0}^{T}f^{-}(s,X_{s},0)\,ds)^{p}\}+\alpha
EK_{\tau_{k}}^{p}.
\end{align}
Taking $\alpha=(2c'(p,\lambda,T))^{-1}$ and combining
(\ref{eq2.2}) with (\ref{eq2.3})  we obtain
\[
E(\int_{0}^{\tau_{k}}|Z_{s}|^{2}\,ds)^{p/2} \le
C(p,\lambda,T)E\{|Y^{*}_{T}|^{p}+(\int_{0}^{T} f_{s}\,ds)^{p}
+(\int_{0}^{T}f^{-}(s,X_{s},0)\,ds)^{p}\}.
\]
Applying Fatou's lemma we conclude from the above inequality and
(\ref{eq2.2}) that
\begin{align*}
E(\int_{0}^{T}|Z_{s}|^{2}\,ds)^{p/2}+EK_{T}^{p} \le
CE\{|Y^{*}_{T}|^{p}+(\int_{0}^{T} f_{s}\,ds)^{p}
+(\int_{0}^{T}f^{-}(s,X_{s},0)\,ds)^{p}\},
\end{align*}
which is the desired estimate.
\end{dow}

\begin{uw}\label{uw1}
Observe that if $f$ does not depend on $z$ then the constant $C$
of Lemma \ref{lm2} depends only on $p$. This follows from the fact
that in this case $c'$ in the key inequality (\ref{eq2.2}) depends
only on $p$.
\end{uw}

\begin{stw}\label{stw1}
Assume that \mbox{\rm(A)} is satisfied and
\[
E(\int_{0}^{T}f^{-}(s,X_{s},0)\,ds)^{p}
+E(\int_{0}^{T}|f(s,0,0)|\,ds)^{p}<+\infty
\]
for some $p>1$ and $X^{+}\in \mathcal{S}^{p}$ such that $X_{t}\ge
Y_{t}$, $t\in [0,T]$. Then if $(Y,Z)$ is a solution of \mbox{\rm
BSDE}$(\xi,f+dK)$ such that $Y\in\mathcal{S}^{p}$, then there
exists $C$ depending only on $\lambda,p,T$ such that for every
$a\ge\mu+\lambda^{2}/[1\wedge(p-1)]$ and every stopping time
$\tau\le T$,
\begin{align*}
&E\sup_{t\le \tau}e^{apt}|Y_{t}|^{p}+E(\int_{0}^{\tau}
e^{2as}|Z_{s}|^{2}\,ds)^{p/2}
+E(\int_{0}^{\tau}e^{as}\,dK_{s})^{p}\\
&\qquad\le CE\bigg(e^{ap\tau}|Y_{\tau}|^{p}
+(\int_{0}^{\tau}e^{as}|f(s,0,0)|\,ds)^{p}+\sup_{t\le \tau}
|e^{at}X^{+}_{t}|^{p}\\
&\qquad\quad+(\int_{0}^{\tau}e^{as}f^{-}(s,X_{s},0)\,ds)^{p}
+(\int_{0}^{\tau}a^{+}e^{as}X^{+}_{s}\,ds)^{p}\bigg).
\end{align*}
Assume additionally that $f$ does not depend on $z$. If $p=1$ and
$X^{+},Y$ are of class \mbox{\rm(D)} then for every $a\ge\mu$,
\begin{align*}
\|e^{a\cdot} Y\|_{1} +E\int_{0}^{T}e^{as}\,dK_{s}&\le
E\bigg(e^{aT}|\xi|
+\int_{0}^{T}e^{as}|f(s,0)|\,ds\\
&\quad+ \int_{0}^{T}e^{as}f^{-}(s,X_{s})\,ds
+\int_{0}^{T}a^{+}e^{as}X^{+}_{s}\,ds\bigg)+\|e^{a\cdot}X^{+}\|_{1}\,.
\end{align*}
\end{stw}
\begin{dow}
To shorten notation we prove the proposition in the case where
$\tau=T$. The proof of the general case requires only minor
technical changes. Moreover, by the change of variables used at
the beginning of the proof of Lemma \ref{lm1} we can reduce the
proof to the case where $a=0$ and
$\mu+\lambda^{2}/[1\wedge(p-1)]\le 0$. Therefore we will assume
that $a,\mu,\lambda$ satisfy the last two conditions.

By It\^o's formula (see Proposition \ref{prop.ito}),
\begin{align*}
|Y_{t}|^{p}&+c(p)\int_{t}^{T}|Y_{s}|^{p-2}\mathbf{1}_{\{Y_{s} \neq
0\}}|Z_{s}|^{2}\,ds= |\xi|^{p}
+p\intt|Y_{s}|^{p-1}\hat{Y_{s}}f(s,Y_{s},Z_{s})\,ds\\
&+p\intt|Y_{s}|^{p-1}\hat{Y}_{s}\,dK_{s}
-p\intt|Y_{s}|^{p-1}\hat{Y_{s}}Z_{s}\,dB_{s},\quad
t\in [0,T].
\end{align*}
By the same method as in the proof of Eq. (6) in \cite{BDHPS} we
deduce from the above inequality that
\begin{align}\label{eq2.4}
\nonumber|Y_{t}|^{p}+\frac{c(p)}{2}
\int_{t}^{T}|Y_{s}|^{p-2}\mathbf{1}_{\{Y_{s}\neq 0\}}|Z_{s}|^{2}\,ds
&\le H-p\intt|Y_{s}|^{p-1}\hat{Y_{s}}Z_{s}\,dB_{s}\\
&\quad+p\intt|Y_{s}|^{p-1}\hat{Y}_{s}\,dK_{s},\quad t\in [0,T],
\end{align}
where $H=|\xi|^{p}+\int_{0}^{T}|Y_{s}|^{p-1}f_{s}\,ds$. Since the
mapping  $\mathbb{R}\ni y\mapsto|y|^{p-1}\hat{y}$ is increasing,
\[
\intt|Y_{s}|^{p-1}\hat{Y}_{s}\,dK_{s} \le
\intt|X^{+}_{s}|^{p-1}\hat X^{+}_{s}\,dK_{s},\quad t\in [0,T].
\]
From this and (\ref{eq2.4}),
\begin{equation}
\label{eq3.06}
|Y_{t}|^{p}+\frac{c(p)}{2}\int_{t}^{T}|Y_{s}|^{p-2}\mathbf{1}_{\{Y_{s}
\neq 0\}}|Z_{s}|^{2}\,ds\le H'
-p\intt|Y_{s}|^{p-1}\hat{Y_{s}}Z_{s}\,dB_{s},
\end{equation}
where $H'=|\xi|^{p}+p\int_{0}^{T}|Y_{s}|^{p-1}f_{s}\,ds
+p\int_{0}^{T}|X^{+}_{s}|^{p-1}\,dK_{s}$. As in the proof of
\cite[Proposition 3.2]{BDHPS} (see (7) and the second inequality
following (8) in \cite{BDHPS}), using the Burkholder-Davis-Gundy
inequality we conclude from (\ref{eq3.06}) that
\begin{equation}\label{eq2.5}
E|Y^{*}_{T}|^{p}\le d_{p} EH'.
\end{equation}
Applying Young's inequality we get
\begin{align}\label{eq2.6}
pd_{p}E\int_{0}^{T}|Y_{s}|^{p-1}f_{t}\,dt \le
pd_{p}E(|Y^{*}_{T}|^{p-1}\int_{0}^{T}f_{t}\,dt) \le
\frac14E|Y^{*}_{T}|^{p}+d'_{p}E(\int_{0}^{T}f_{t}\,dt)^{p}
\end{align}
and
\begin{align}\label{eq2.7}
pd_{p}E\int_{0}^{T}|X^{+}_{t}|^{p-1}\,dK_{t} \le
d'(p,\alpha)E|X^{+,*}_{T}|^{p}+\alpha EK^{p}_{T}\,.
\end{align}
By Lemma \ref{lm1}, there exists $d(p,\lambda,T)>0$ such that
\[
EK_{T}^{p}\le d(p,\lambda,T)E\{|Y^{*}_{T}|^{p}
+(\int_{0}^{T}|Z_{s}|^{2}\,ds)^{p/2}
+(\int_{0}^{T}f^{-}(s,X_{s},0)\,ds)^{p}\}.
\]
From this  and Lemma \ref{lm2} we see that there exists
$c(p,\lambda,T)>0$ such that
\begin{align}
\label{eq2.8} EK^{p}_{T}\le c(p,\lambda,T)E\{|Y^{*}_{T}|^{p}
+(\int_{0}^{T}f_{s}\,ds)^{p}
+(\int_{0}^{T}f^{-}(s,X_{s},0)\,ds)^{p}\}.
\end{align}
Put $\alpha=(4c(p,\lambda,T))^{-1}$. Then from
(\ref{eq2.5})--(\ref{eq2.8}) it follows that there is
$C(p,\lambda,T)$ such that
\begin{align*}
E|Y_{T}^{*}|^{p}&\le C(p,\lambda,T)E\{|\xi|^{p}
+(\int_{0}^{T}|f(s,0,0)|\,ds)^{p} \\
&\quad+ (\int_{0}^{T}f^{-}(s,X_{s},0)\,ds)^{p}+\sup_{t\le T}
|X^{+}_{t}|^{p}\}.
\end{align*}
Hence, by (\ref{eq2.8}) and  Lemma \ref{lm2},
\begin{align*}
E|Y^{*}_{\tau}|^{p}+E(\int_{0}^{\tau}|Z_{s}|^{2}\,ds)^{p/2}
+EK_{T}^{p}&\le CE\bigg(|Y_{\tau}|^{p}
+(\int_{0}^{\tau}|f(s,0,0)|\,ds)^{p}
+|X^{+,*}_{\tau}|^{p}\\
&\quad+(\int_{0}^{\tau}f^{-}(s,X_{s},0)\,ds)^{p}\bigg).
\end{align*}
From this  the first assertion follows. Now suppose that $f$ does
not depend on $z$. As in the first part of the proof we may assume
that $\mu\le 0$ and $a=0$. Applying It\^o's formula (see
Proposition \ref{prop.ito}) we conclude that for any stopping
times $\sigma\le\tau\le T$,
\begin{align}\label{eq2.9}
|Y_{\sigma}|\le|Y_{\tau}|+\int_{\sigma}^{\tau}f(s,Y_{s})\hat{Y}_{s}\,ds
+\int_{\sigma}^{\tau}\hat{Y}_{s}\,dK_{s}
-\int_{\sigma}^{\tau}Z_{s}\hat{Y}_{s}\,dB_{s}.
\end{align}
Let us define $\tau_k$ by (\ref{eq3.01}). Then
$\int_{0}^{\tau_{k}\wedge\cdot}Z_{s}\hat{Y}_{s}\,dB_{s}$ is a
uniformly integrable martingale. Using this, the fact that
$Y$ is of class (D) and monotonicity of $f$ with respect to $y$ we deduce
from (\ref{eq2.9}) that $|Y_{\sigma}|\le
E(|\xi|+\int_{0}^{T}|f(s,0)|\,ds +K_{T}|\mathcal{F}_{\sigma})$,
hence that
\begin{equation}
\label{eq2.10}
\|Y\|_{1}\le E(|\xi|+\int_{0}^{T}|f(s,0)|\,ds+K_{T}).
\end{equation}
On the other hand, $-f(t,Y_{t})\le -f(t,X_{t})$ for $t\in [0,T]$
since $Y_{t}\le X_{t}$, $t\in [0,T]$. Therefore
\begin{align*}
K_{\tau}&=Y_{0}-Y_{\tau}-\int_{0}^{\tau}f(s,Y_{s})\,ds
+\int_{0}^{\tau}Z_{s}\,dB_{s}\\
&\le X_{0}-Y_{\tau} -\int_{0}^{\tau}f(s,X_{s})\,ds
+\int_{0}^{\tau}Z_{s}\,dB_{s}.
\end{align*}
Taking $\tau=\tau_{k}$ and using the fact that $Y$ is of class (D)
we deduce from the above inequality that
\[
EK_{T}\le EX^{+}_{0}+E|\xi|+E\int_{0}^{T}f^{-}(s,X_{s})\,ds.
\]
Combining this with (\ref{eq2.10}) we get the desired result.
\end{dow}

\begin{uw}
If $f$ does not depend on $z$ then the constant $C$ of the first
assertion of Proposition \ref{stw1} depends only on $p$. To see
this it suffices to observe that if $f$ does not depend on $z$
then the constant $c$ in the key inequality (\ref{eq2.8}) depends
only on $p$ (see Remark \ref{uw1}).
\end{uw}

\nsubsection{Some useful tools}
\label{sec4}

We begin with a useful comparison result for solutions of
(\ref{eq3.02}) with $K\equiv0$.

\begin{stw}\label{stw2}
Let $(Y^{1},Z^{1}), (Y^{2},Z^{2})$ be solutions of \mbox{\rm
BSDE}$(\xi^{1},f^{1})$,  \mbox{\rm BSDE}$(\xi^{2},f^{2})$,
respectively. Assume that $(Y^{1}-Y^{2})^{+}\in\mathcal{S}^{q}$
for some $q>1$. If $\xi^{1}\le \xi^{2}$ and for a.e. $t\in [0,T]$
either
\begin{equation}
\label{eq4.01} \mathbf{1}_{\{Y^{1}_{t}>Y^{2}_{t}\}}
(f^{1}(t,Y^{1}_{t},Z^{1}_{t})-f^{2}(t,Y^{1}_{t},Z^{1}_{t}))\le
0,\quad f^{2}\mbox{ satisfies }\mbox{\rm(H2)}, \mbox{\rm(H3)}
\end{equation}
or
\begin{equation}
\label{eq4.02} \mathbf{1}_{\{Y^{1}_{t}>Y^{2}_{t}\}}
(f^{1}(t,Y^{2}_{t},Z^{2}_{t})-f^{2}(t,Y^{2}_{t},Z^{2}_{t}))\le
0, \quad f^{1}\mbox{ satisfies }\mbox{\rm(H2)}, \mbox{\rm(H3)}
\end{equation}
is satisfied then $Y^{1}_{t}\le Y^{2}_{t}$, $t\in [0,T]$.
\end{stw}
\begin{dow}
We show the proposition  in case (\ref{eq4.01}) is satisfied. If
 (\ref{eq4.02}) is satisfied, the proof is analogous. Without
loss of generality we may assume that $\mu\le 0$. By the
It\^o-Tanaka formula,  for every $p\in(1,q)$ and every stopping
time $\tau\le T$,
\begin{align}
\label{eq3.1} &|(Y^{1}_{t\wedge\tau}-Y^{2}_{t\wedge\tau})^{+}|^{p}
+\frac{p(p-1)}{2}\int_{t\wedge\tau}^{\tau}\mathbf{1}_{\{Y^{1}_{s}
\neq Y^{2}_{s}\}}|(Y^{1}_{s}-Y^{2}_{s})^{+}|^{p-2}
|Z^{1}_{s}-Z^{2}_{s}|^{2}\,ds\nonumber\\
&\qquad=|(Y^{1}_{\tau}-Y^{2}_{\tau})^{+}|^{p}
+p\int_{t\wedge\tau}^{\tau}|(Y^{1}_{s}-Y^{2}_{s})^{+}|^{p-1}
(f^{1}(s,Y^{1}_{s},Z^{1}_{s})-f^{2}(s,Y^{2}_{s},Z^{2}_{s}))\,ds
\nonumber\\
&\qquad\quad-p\int_{t\wedge\tau}^{\tau}|(Y^{1}_{s}-Y^{2}_{s})^{+}|^{p-1}
(Z^{1}_{s}-Z^{2}_{s})\,dB_{s}.
\end{align}
By (\ref{eq4.01}),
\begin{align*}
&\mathbf{1}_{\{Y^{1}_{t}>Y^{2}_{t}\}}
(f^{1}(t,Y^{1}_{t},Z^{1}_{t})-f^{2}(t,Y^{2}_{t},Z^{2}_{t}))\\
&\qquad=\mathbf{1}_{\{Y^{1}_{t}>Y^{2}_{t}\}}
(f^{1}(t,Y^{1}_{t},Z^{1}_{t})-f^{2}(t,Y^{1}_{t},Z^{1}_{t}))\\
&\qquad\quad+\mathbf{1}_{\{Y^{1}_{t}>Y^{2}_{t}\}}
(f^{2}(t,Y^{1}_{t},Z^{1}_{t})
-f^{2}(t,Y^{2}_{t},Z^{2}_{t}))\\
&\qquad\le\mathbf{1}_{\{Y^{1}_{t}>Y^{2}_{t}\}}
(f^{2}(t,Y^{1}_{t},Z^{1}_{t})-f^{2}(t,Y^{2}_{t},Z^{1}_{t}))\\
&\qquad\quad+\mathbf{1}_{\{Y^{1}_{t}>Y^{2}_{t}\}}
(f^{2}(t,Y^{2}_{t},Z^{1}_{t})-f^{2}(t,Y^{2}_{t},Z^{2}_{t}))\\
&\qquad\le\lambda\mathbf{1}_{\{Y^{1}_{t}>Y^{2}_{t}\}}
|Z^{1}_{t}-Z^{2}_{t}|.
\end{align*}
From this, (\ref{eq3.1}) and Young's inequality,
\begin{align*}
&|(Y^{1}_{t\wedge\tau}-Y^{2}_{t\wedge\tau})^{+}|^{p}+\frac{p(p-1)}{2}
\int_{t\wedge\tau}^{\tau}\mathbf{1}_{\{Y^{1}_{s}\neq Y^{2}_{s}\}}|(Y^{1}_{s}-Y^{2}_{s})^{+}|^{p-2}
|Z^{1}_{s}-Z^{2}_{s}|^{2}\,ds\\
&\qquad\le |(Y^{1}_{\tau}-Y^{2}_{\tau})^{+}|^{p}
+p\lambda\int_{t\wedge\tau}^{\tau}|(Y^{1}_{s}
-Y^{2}_{s})^{+}|^{p-1}|Z^{1}_{s}-Z^{2}_{s}|\,ds\\
&\qquad\quad-p\int_{t\wedge\tau}^{\tau}|(Y^{1}_{s}-Y^{2}_{s})^{+}|^{p-1}
(Z^{1}_{s}-Z^{2}_{s})\,dB_{s}\\
&\qquad\le|(Y^{1}_{\tau}-Y^{2}_{\tau})^{+}|^{p}
+\frac{p\lambda^{2}}{p-1}\int_{t\wedge\tau}^{\tau}
|(Y^{1}_{s}-Y^{2}_{s})^{+}|^{p}\,ds\\
&\qquad\quad+\frac{p(p-1)}{4}\int_{t\wedge\tau}^{\tau}\mathbf{1}_{\{Y^{1}_{s}
\neq Y^{2}_{s}\}}
|(Y^{1}_{s}-Y^{2}_{s})^{+}|^{p-2}|Z^{1}_{s}-Z^{2}_{s}|^{2}\,ds
\\&\qquad\quad-p\int_{t\wedge\tau}^{\tau}|(Y^{1}_{s}-Y^{2}_{s})^{+}|^{p-1}
(Z^{1}_{s}-Z^{2}_{s})\,dB_{s}.
\end{align*}
Let $\tau_{k}=\inf\{t\in [0,T];\int_{0}^{t}
|(Y^{1}_{s}-Y^{2}_{s})^{+}|^{2(p-1)}
|Z^{1}_{s}-Z^{2}_{s}|^{2}\,ds\ge k\}\wedge T$. From the above
estimate it follows that
\[
E|(Y^{1}_{t\wedge\tau_{k}}-Y^{2}_{t\wedge\tau_{k}})^{+}|^{p} \le
E|(Y^{1}_{\tau_{k}}-Y^{2}_{\tau_{k}})^{+}|^{p}
+\frac{p\lambda^{2}}{p-1}E\int_{t\wedge\tau_{k}}^{\tau_{k}}
|(Y^{1}_{s}-Y^{2}_{s})^{+}|^{p}\,ds.
\]
Since$(Y^{1}-Y^{2})^{+}\in \mathcal{S}^{q}$, letting
$k\rightarrow\infty$ and using the assumptions we get
\[
E|(Y^{1}_{t}-Y^{2}_{t})^{+}|^{p} \le
\frac{p\lambda^{2}}{p-1}E\int_{t}^{T}
|(Y^{1}_{s}-Y^{2}_{s})^{+}|^{p}\,ds,\quad t\in [0,T].
\]
By Gronwall's lemma, $E|(Y^{1}_{t}-Y^{2}_{t})^{+}|^{p}=0,\, t\in [0,T]$, from
which the desired result follows.
\end{dow}

\begin{lm}\label{lm3}
Assume that $\{(X^{n},Y^{n},K^{n})\}$ is a sequence of  real
valued c\`adl\`ag progressively measurable processes such that
\begin{enumerate}
\item [\rm (a)]$Y^{n}_{t}=-K^{n}_{t}
 +X_{t}^{n},\, t\in [0,T],\,K^{n}$-increasing, $K^{n}_{0}=0,$
\item [\rm (b)] $Y^{n}_{t}\uparrow Y_{t}$, $t\in [0,T]$,
$Y^{1},Y$ are of class (D),
\item [\rm (c)]There exists a c\`adl\`ag process $X$ such that for some
subsequence $\{n'\}$, $X^{n'}_{\tau}\rightarrow X_{\tau}$ weakly
in $\mathbb{L}^{1}(\FF_{T})$ for every stopping time $\tau\le T$.
\end{enumerate}
Then $Y$ is c\`adl\`ag and there exists a c\`adl\`ag increasing
process $K$ such that $K^{n'}_{\tau}\rightarrow K_{\tau}$ weakly
in $\mathbb{L}^{1}(\FF_{T})$ for every stopping time $\tau\le T$ and
\[
Y_{t}=-K_{t}+X_{t},\quad t\in [0,T].
\]
\end{lm}
\begin{dow}
From (b) it follows  that $Y^{n'}_{\tau}\rightarrow Y_{\tau}$
weakly  in $\mathbb{L}^{1}(\FF_{T})$ for every stopping time
$\tau\le T$. Set $K_{t}=X_{t}-Y_{t}$. By the above and
(c), $K^{n'}_{\tau}\rightarrow K_{\tau}$ weakly in
$\mathbb{L}^{1}(\Omega)$ for every stopping time $\tau\le T$. If
$\sigma,\tau$ are stopping times such that $\sigma\le\tau\le T$
then $K_{\sigma}\le K_{\tau}$ since $K^{n}_{\sigma}\le
K^{n}_{\tau}$, $n\in\BN$. Therefore $K$ is increasing. The fact
that $Y,K$ are c\`adl\`ag processes follows easily from
\cite[Lemma 2.2]{Peng}.
\end{dow}
\medskip

In what follows we say that a sequence $\{\tau_{k}\}$ of stopping
times is stationary if
\[
P(\liminf_{k\rightarrow+\infty} \{\tau_{k}=T\})=1.
\]

\begin{lm}\label{lm4}
Assume that $\{Y^{n}\}$ is a nondecreasing sequence of  continuous
processes such that $\sup_{n\ge1}E|Y^{n,*}_{T}|^{q}<\infty$ for
some $q>0$. Then there exists a stationary sequence $\{\tau_{k}\}$
of stopping times such that $Y^{n,*}_{\tau_{k}}\le
k\vee|Y^{n}_{0}|$, $P$-a.s. for every $k\in\mathbb{N}$.
\end{lm}
\begin{dow}
Set $V^{n}_{t}=\sup_{0\le s\le t}(Y^{n}_{s}-Y^{1}_{s})$. Then
$V^{n}$ is nonnegative and $V^{n}\in\mathcal{V}^{+}_{c}$. Since
$\{Y^{n}\}$ is nondecreasing, there exists an increasing process
$V$ such that $V^{n}_{t}\uparrow V_{t}$, $t\in [0,T]$. By Fatou's
lemma,
\[
EV^{q}_{T}\le \liminf_{n\rightarrow+\infty} E|V^{n}_{T}|^{q}\le
c(q)\sup_{n\ge1}E|Y^{n,*}_{T}|^{q}<\infty.
\]
Now, set $V'_{t}=\inf_{t<t'\le T}V_{t'}$, $t\in [0,T]$ and then
$\tau_{k}=\inf\{t\in [0,T]; Y^{1,*}_{t}+V'_{t}>k\}\wedge T$. It is
known that $V'$ is a progressively measurable c\`adl\`ag process.
Since $V_{T}$ is integrable, the sequence $\{\tau_{k}\}$ is
stationary. From the above it follows that if $\tau_{k}>0$ then
\[
Y^{n,*}_{\tau_{k}}=Y^{n,*}_{\tau_{k}-} \le
V'_{\tau_{k}-}+Y^{1,*}_{\tau_{k}-}\le k,\quad k\in\mathbb{N},
\]
and the proof is complete.
\end{dow}

\begin{lm}\label{lm5}
If $\{Z^{n}\}$ is a sequence of progressively measurable processes
such that
$\sup_{n\ge1}E(\int_{0}^{T}|Z^{n}_{t}|^{2}\,dt)^{p/2}<\infty$ for
some $p>1$, then there exists $Z\in M^{p}$ and a subsequence
$\{n'\}$ such that for every stopping time $\tau\le T$,
$\int_{0}^{\tau}Z^{n'}_{t}\,dB_{t}\rightarrow
\int_{0}^{\tau}Z_{t}\,dB_{t}$ weakly in $\mathbb{L}^{p}(\FF_{T})$.
\end{lm}
\begin{dow}
Since $\{Z^{n}\}$ is bounded in $\mathbb{L}^{2,p}(\FF)$ and the
space $\mathbb{L}^{2,p}(\FF)$ is reflexive, there exists a
subsequence (still denoted by $\{n\}$) and
$Z\in\mathbb{L}^{2,p}(\FF)$ such that $Z^{n}\rightarrow Z$ weakly
in $\mathbb{L}^{2,p}(\FF)$. It is known that if
$\xi\in\mathbb{L}^{p'}(\FF_{T})$, where $p'=p/(p-1)$, then there
exists $\eta\in
\mathbb{L}^{2,p'}(\FF)=(\mathbb{L}^{2,p}(\FF))^{*}$ such that
\begin{equation}
\label{eq3.2}
\xi=E\xi+\int_{0}^{T}\eta_{t}\,dB_{t}.
\end{equation}
Let $f\in (\mathbb{L}^{p}(\FF_{T}))^{*}$. Then there exists
$\xi\in\mathbb{L}^{p'}(\FF_{T})$ such that $f(\zeta)=E\zeta\xi$
for every $\zeta\in\mathbb{L}^{p}(\FF_{T})$. Let
$\eta\in\mathbb{L}^{2,p'}(\FF)$ be such that (\ref{eq3.2}) is
satisfied. Without loss of generality we may assume that $E\xi=0$.
Then by  It\^o's isometry,
\begin{align*}
f(\int_{0}^{T}Z^{n}_{t}\,dB_{t}) &=E\xi\int_{0}^{T}
Z^{n}_{t}\,dB_{t}
=E\int_{0}^{T}\eta_{t}\,dB_{t}\int_{0}^{T}Z^{n}_{t}\,dB_{t}\\
%=E\langle \int_{0}^{\cdot}\eta_{s}\,dB_{s},
%\int_{0}^{\cdot}Z^{n}_{s}\,dB_{s}\rangle_{T}\\&
&=E\int_{0}^{T}\eta_{t}Z^{n}_{t}\,dt\rightarrow
E\int_{0}^{T}\eta_{t}Z_{t}=f(\int_{0}^{T}Z_{t}\,dB_{t}).
\end{align*}
Since the same reasoning applies to the sequence
$\{\mathbf{1}_{\{\cdot\le\tau\}}Z^{n}\}$ in place of $\{Z^n\}$,
the lemma follows.
\end{dow}

\nsubsection{Existence and uniqueness results for $p>1$}
\label{sec5}

First we recall the definition of a solution $(Y,Z,K)$ of
(\ref{eq1.1}). Note that a priori we do not impose any
integrability conditions on the processes $Y,Z,K$.

\begin{df}
We say that a triple $(Y,Z,K)$ of progressively measurable
processes  is a solution of RBSDE$(\xi,f,L)$ iff
\begin{enumerate}
\item [\rm(a)]$K$ is an increasing continuous process, $K_{0}=0$,
\item [\rm(b)]$Z\in M$ and the mapping $[0,T]\ni t\mapsto
f(t,Y_{t},Z_{t})$ belongs to $\mathbb{L}^{1}(0,T),\, P$-a.s.,
\item [\rm(c)]$Y_{t}=\xi+\int_{t}^{T}f(s,Y_{s},Z_{s})\,ds
+\int_{t}^{T}dK_{s}-\int_{t}^{T}Z_{s}\,dB_{s},\quad t\in [0,T],$
\item [\rm(d)]$Y_{t}\ge L_{t},\, t\in [0,T]$,
$\int_{0}^{T}(Y_{t}-L_{t})\,dK_{t}=0.$
\end{enumerate}
\end{df}

\begin{stw}\label{stw2.5}
Let $(Y^{1},Z^{1},K^{1}), (Y^{2},Z^{2},K^{2})$ be solutions of
\mbox{\rm RBSDE}$(\xi^{1},f^{1},L^{1})$, \mbox{\rm
RBSDE}$(\xi^{2},f^{2},L^{2})$, respectively. Assume that
$(Y^{1}-Y^{2})^{+}\in\mathcal{S}^{q}$ for some $q>1$. If
$\xi^{1}\le\xi^{2}$, $L^{1}_{t}\le L^{2}_{t}$, $t\in [0,T]$, and
either \mbox{\rm (\ref{eq4.01})} or \mbox{\rm(\ref{eq4.02})} is
satisfied then $Y^{1}_{t}\le Y^{2}_{t}$, $t\in [0,T]$.
\end{stw}
\begin{dow}
Assume that (\ref{eq4.01}) is satisfied. Let $q>1$ be such that
$(Y^{1}-Y^{2})^{+}\in\mathcal{S}^{q}$. Without loss of generality
we may assume that $\mu\le 0$. By the It\^o-Tanaka formula, for
$p\in(1,q)$ and every stopping time $\tau\le T$,
\begin{align}\label{eqt}
\nonumber &|(Y^{1}_{t\wedge\tau}-Y^{2}_{t\wedge\tau})^{+}|^{p}
+\frac{p(p-1)}{2}\int_{t\wedge\tau}^{\tau}\mathbf{1}_{\{Y^{1}
\neq Y^{2}_{s}\}}|(Y^{1}_{s}-Y^{2}_{s})^{+}|^{p-2}
|Z^{1}_{s}-Z^{2}_{s}|^{2}\,ds\nonumber\\
&\quad=|(Y^{1}_{\tau}-Y^{2}_{\tau})^{+}|^{p}
+p\int_{t\wedge\tau}^{\tau}|(Y^{1}_{s}-Y^{2}_{s})^{+}|^{p-1}
(f^{1}(s,Y^{1}_{s},Z^{1}_{s})-f^{2}(s,Y^{2}_{s},Z^{2}_{s}))\,ds\nonumber\\
&\qquad+p\int_{t\wedge\tau}^{\tau} |(Y^{1}_{s}-Y^{2}_{s})^{+}|^{p-1 }\,
(dK^{1}_{s}-dK^{2}_{s})\nonumber\\
&\qquad-p\int_{t\wedge\tau}^{\tau}|(Y^{1}_{s}-Y^{2}_{s})^{+}|^{p-1}
(Z^{1}_{s}-Z^{2}_{s})\,dB_{s}.
\end{align}
By  monotonicity of the function $x\mapsto \hat{x}|x|^{p-1}$,
condition (d)  of the definition of a solution of reflected BSDE
and the fact that $L^1_t\le L^2_t$ for $t\in[0,T]$,
\begin{align*}
\int_{t\wedge\tau}^{\tau}|(Y^{1}_{s}-Y^{2}_{s})^{+}|^{p-1}\,
(dK^{1}_{s}-dK^{2}_{s})&\le \int_{t\wedge\tau}^{\tau}
|(Y^{1}_{s}-Y^{2}_{s})^{+}|^{p-1}\,dK^{1}_{s}\\&
\le \int_{t\wedge\tau}^{\tau} |(Y^{1}_{s}-L^{1}_{s})^{+}|^{p-1}\,dK^{1}_{s}=0.
\end{align*}
Combining this with (\ref{eqt}) we get estimate (\ref{eq3.1}) in
Proposition \ref{stw2}. Therefore repeating arguments following
(\ref{eq3.1}) in the proof of that proposition we obtain the
desired result. The proof in case (\ref{eq4.02}) is satisfied is
analogous and therefore left to the reader.
\end{dow}

\begin{stw}\label{stw3}
If $f$ satisfies \mbox{\rm(H2)}, \mbox{\rm(H3)} then there exists
at most one solution $(Y,Z,K)$ of \mbox{\rm RBSDE}$(\xi,f,L)$ such
that $Y\in \mathcal{S}^{p}$ for some $p>1$.
\end{stw}
\begin{dow}
Follows immediately from Proposition \ref{stw2.5} and uniqueness
of the Doob-Meyer decomposition of semimartingales.
\end{dow}

\begin{tw}\label{tw1}
Let $p>1$.
\begin{enumerate}
\item[\rm(i)] Assume  \mbox{\rm(H1)}--\mbox{\rm(H6)}.
Then there exists a solution $(Y,Z,K)$ of \mbox{\rm
RBSDE}$(\xi,f,L)$ such that $(Y,Z,K)\in \mathcal{S}^{p}\otimes
M^{p}\otimes\mathcal{V}^{+,p}_{c}$ iff \mbox{\rm(H7)} is
satisfied.
\item[\rm(ii)]Assume \mbox{\rm(H1)}--\mbox{\rm(H7)}. For
$n\in\mathbb{N}$ let $(Y^{n},Z^{n})$ be a solution of the BSDE
\begin{equation}
\label{eq5.01}
Y^{n}_{t}=\xi+\int_{t}^{T}f(s,Y^{n}_{s},Z^{n}_{s})\,ds
+\int_{t}^{T}dK^n_s-\int_{t}^{T}Z^{n}_{s}\,dB_{s},\, t\in [0,T]
\end{equation}
with
\begin{equation}
\label{eq5.02} K_{t}^{n}=\int_{0}^{t}n(Y^{n}_{s}-L_{s})^{-}\,ds
\end{equation}
such that $(Y^{n},Z^{n})\in\mathcal{S}^{p}\otimes M^{p}$.
Then
\begin{equation}
\label{eq5.03} E\sup_{t\le T}|Y^{n}_{t}-Y_{t}|^{p} +E\sup_{t\le
T}|K^{n}_{t}-K_{t}|^{p}
+E(\int_{0}^{T}|Z^{n}_{t}-Z_{t}|^{2}dt)^{p/2}\rightarrow 0
\end{equation}
as $n\rightarrow +\infty$.
\end{enumerate}
\end{tw}
\begin{dow}
Without loss of generality we may assume that $\mu\le0$. Assume
that there is a solution $(Y,Z,K)\in\mathcal{S}^{p}\otimes
M^{p}\otimes\mathcal{V}^{+,p}_{c} $ of RBSDE$(\xi,f,L)$. Then by
\cite[Remark 4.3]{BDHPS},
\begin{equation}\label{eq4.1}
E(\int_{0}^{T}|f(s,Y_{s},Z_{s})|\,ds)^{p}\le cE\{|\xi|^{p}
+(\int_{0}^{T}f_{s}\,ds)^{p}+K^{p}_{T}\}
\end{equation}
which in view of (H2) and the fact that $Y_{t}\ge L_{t}$, $t\in
[0,T]$ shows (H7). Conversely, let us assume that (H1)--(H7) are
satisfied. Let $(Y^{n},Z^{n})$ be a solution of (\ref{eq5.01}) such that $(Y^{n},Z^{n})\in
\mathcal{S}^{p}\otimes M^{p}$. We will show that there exists a
process $\overline{X}\in\mathcal{H}^{p}_{c}$ such that
$\overline{X}_{t}\ge Y^{n}_{t}$, $t\in [0,T]$ for every
$n\in\mathbb{N}$. Since $X\in\HH^{p}_{c}$, there exist $M\in
\MM^{p}_{c}$ and $V\in\mathcal{V}^{p}_{c}$  such that $X=V+M$. By
the representation property of Brownian filtration, there exists
$Z'\in M^{p}$ such that
\[
X_{t}=X_{T}-\intt dV_{s}-\intt Z'_{s}\,dB_{s},\quad t\in [0,T].
\]
The above identity can be rewritten in the form
\begin{align*}
X_{t}&=X_{T}+\intt f(s,X_{s},Z'_{s})\,ds
-\intt (f^{+}(s,X_{s},Z'_{s})\,ds+dV^{+}_{s})\\&\quad
+\intt (f^{-}(s,X_{s},Z'_{s})\,ds+dV^{-}_{s})
-\intt Z'_{s}\,dB_{s},\quad t\in [0,T].
\end{align*}
By \cite[Theorem 4.2]{BDHPS}, there exists a unique solution
$(\overline{X},\overline{Z})\in\mathcal{S}^{p}\otimes M^{p}$ of
the BSDE
\begin{align*}
\overline{X}_{t}=\xi\vee X_T +\intt
f(s,\overline{X}_{s},\overline{Z}_{s})\,ds +\intt
(f^{-}(s,X_{s},Z'_{s})\,ds+dV^{-}_{s})-\intt
\overline{Z}_{s}\,dB_{s}.
\end{align*}
By Proposition \ref{stw2}, $\overline{X}_{t}\ge X_{t}\ge L_{t},\,
t\in [0,T]$. Hence
\begin{align*}
\overline{X}_{t}&=\xi\vee X_T +\intt
f(s,\overline{X}_{s},\overline{Z}_{s})\,ds +\intt
n(\overline{X}_{s}-L_{s})^{-}\,ds
\\&\quad+\intt (f^{-}(s,X_{s},Z'_{s})\,ds
+dV^{-}_{s})-\intt \overline{Z}_{s}\,dB_{s},\quad t\in [0,T],
\end{align*}
so using once again Proposition \ref{stw2} we see that
$\overline{X}_{t}\ge Y^{n}_{t}$, $t\in [0,T]$. By \cite[Remark
4.3]{BDHPS},
$E(\int_{0}^{T}|f(s,\overline{X}_{s},0)|\,ds)^{p}<\infty$. Hence,
by Lemma \ref{lm1} and Proposition \ref{stw1},
\begin{align}
\label{eq4.2}
&E|Y^{n,*}_{T}|^{p}+E(\int_{0}^{T}|Z^{n}_{s}|^{2}\,ds)^{p/2}
+E|K^{n}_{T}|^{p}\nonumber\\
&\qquad\le C(p,\lambda,T)E\{|\xi|^{2} +(\int_{0}^{T}f_{s}\,ds)^{p}
+(\int_{0}^{T}|f(s,\overline{X}_{s},0)|\,ds)^{p}\}.
\end{align}
From this and \cite[Remark 4.3]{BDHPS},
\begin{equation}
\label{T1} E(\int_{0}^{T}|f(s,Y^{n}_{s},Z_{s}^{n})|\,ds)^{p}\le
C^{'}(p,\lambda,T).
\end{equation}
By Proposition \ref{stw2} there exists a progressively measurable
process $Y$ such that $Y^{n}_{t}\uparrow Y_{t}$, $t\in [0,T]$.
Using the monotone  convergence of $Y^{n}$, (H3)--(H5),
(\ref{eq4.2}), (\ref{T1}) and the Lebesgue dominated
convergence theorem we conclude that
\begin{equation}
\label{T2}
E(\int_{0}^{T}|f(s,Y^{n}_{s},0)-f(s,Y_{s},0)|\,ds)^{p}\rightarrow
0
\end{equation}
Moreover, by (H2) and (\ref{eq4.2}),
\[
\sup_{n\ge1}E\int_{0}^{T}
|f(s,Y^{n}_{s},Z^{n}_{s})-f(s,Y^{n}_{s},0)|^{p}\,ds<\infty.
\]
It follows in particular that there exists a  process
$\eta\in\mathbb{L}^{p}(\FF)$ such that
\[
\int_{0}^{\tau}(f(s,Y^{n}_{s},Z^{n}_{s})-f(s,Y^{n}_{s},0))\,ds
\rightarrow \int_{0}^{\tau}\eta_{s}\,ds
\]
weakly in $\mathbb{L}^{1}(\FF_{T})$ for every stopping time
$\tau\le T$. Consequently,  by Lemmas \ref{lm3} and \ref{lm5}, $Y$
is a c\`adl\`ag process and there exist $Z\in M^{p}$ and a
c\`adl\`ag increasing process $K$ such that $K_{0}=0$ and
\begin{equation}
\label{eq4.3}
Y_{t}=\xi+\intt f(s,Y_{s},0)\,ds+\intt\eta_{s}\,ds+\intt dK_{s}
-\intt Z_{s}\,dB_{s},\quad t\in [0,T].
\end{equation}
From (\ref{eq5.01}), (\ref{eq4.2}), (\ref{T1}) and pointwise
convergence of the sequence $\{Y^{n}\}$ one can deduce that
$E\int_{0}^{T}(Y_{s}-L_{s})^{-}\,ds=0$, which when combined with
(H6) and the fact that $Y$ is c\`adl\`ag implies that $Y_{t}\ge
L_{t}$, $t\in [0,T]$. From this, the monotone character of the
convergence of the sequence $\{Y^{n}\}$ and Dini's theorem we
conclude that
\begin{align}\label{eq4.4}
E|(Y^{n}-L)^{-,*}_{T}|^{p}\rightarrow 0.
\end{align}
By Proposition \ref{prop.ito}, for $n,m\in\mathbb{N}$ we have
\begin{align}\label{eq4.5}
\nonumber &|Y^{n}_{t}-Y^{m}_{t}|^{p}+c(p)\intt
|Y^{n}_{s}-Y^{m}_{s}|^{p-2} \mathbf{1}_{\{Y^{n}_{s}-Y^{m}_{s}\neq
0\}} |Z^{n}_{s}-Z^{m}_{s}|^{2}\,ds\\
\nonumber &\qquad = p\intt |Y^{n}_{s}-Y^{m}_{s}|^{p-1}
\widehat{Y^{n}_{s}-Y^{m}_{s}}
(f(s,Y^{n}_{s},Z^{n}_{s})-f(s,Y^{m}_{s},Z^{m}_{s}))\,ds\\
\nonumber &\qquad\quad+p\intt|Y^{n}_{s}-Y^{m}_{s}|^{p-1}
\widehat{Y^{n}_{s}-Y^{m}_{s}}(dK^{n}_{s}-dK^{m}_{s})\\
&\qquad\quad -p\intt |Y^{n}_{s}-Y^{m}_{s}|^{p-1}
\widehat{Y^{n}_{s}-Y^{m}_{s}}(Z^{n}_{s}-Z^{m}_{s})\,dB_{s},\quad
t\in [0,T].
\end{align}
By monotonicity of the function $\mathbb{R}\ni
x\mapsto|x|^{p-1}\hat{x}$,
\begin{align}
\label{eq4.6}
\intt|Y^{n}_{s}-Y^{m}_{s}|^{p-1}\widehat{Y^{n}_{s}-Y^{m}_{s}}\,dK^{n}_{s}
\le \intt|(Y^{m}_{s}-L_{s})^{-}|^{p-1}
\widehat{(Y^{m}_{s}-L_{s})^{-}}\,dK^{n}_{s}
\end{align}
and
\begin{align}
\label{eq4.7}
-\intt|Y^{n}_{s}-Y^{m}_{s}|^{p-1}
\widehat{Y^{n}_{s}-Y^{m}_{s}}\,dK^{m}_{s}
\le \intt|(Y^{n}_{s}-L_{s})^{-}|^{p-1}
\widehat{(Y^{n}_{s}-L_{s})^{-}}\,dK^{m}_{s}.
\end{align}
By (H2), (H3), (\ref{eq4.5})--(\ref{eq4.7}) and H\"older's
inequality,
\begin{align}
\label{eq.tem2012}
\nonumber & E|Y^{n}_{t}-Y^{m}_{t}|^{p}
+c(p)E\intt|Y^{n}_{s}-Y^{m}_{s}|^{p-2}
\mathbf{1}_{\{Y^{n}_{s}-Y^{m}_{s}\neq0\}}
|Z^{n}_{s}-Z^{m}_{s}|^{2}\,ds\\
&\nonumber\quad\le p\lambda E\int_{0}^{T}
|Y^{n}_{s}-Y^{m}_{s}|^{p-1}|Z^{n}_{s}-Z^{m}_{s}|\,ds
+(E|(Y^{n}-L)^{-,*}_{T}|^{p})^{(p-1)/p}(E|K^{m}_{T}|^{p})^{1/p}
\\&\qquad
+(E|(Y^{m}-L)^{-,*}_{T}|^{p})^{(p-1)/p}(E|K^{n}_{T}|^{p})^{1/p}.
\end{align}
Since
\begin{align*}
p\lambda |Y^{n}_{s}-Y^{m}_{s}|^{p-1}|Z^{n}_{s}-Z^{m}_{s}|
&\le\frac{p\lambda^{2}}{1\wedge(p-1)}|Y^{n}_{s}-Y^{m}_{s}|^{p}\\&
\quad+\frac{c(p)}{2}\mathbf{1}_{\{Y^{n}_{s}-Y^{m}_{s} \neq
0\}}|Y^{n}_{s}-Y^{m}_{s}|^{p-2}|Z^{n}_{s}-Z^{m}_{s}|^{2},
\end{align*}
from (\ref{eq.tem2012}) we get
\begin{align*}
\label{eq.tem2012}
\nonumber & E|Y^{n}_{t}-Y^{m}_{t}|^{p}
+\frac{c(p)}{2}E\intt|Y^{n}_{s}-Y^{m}_{s}|^{p-2}
\mathbf{1}_{\{Y^{n}_{s}-Y^{m}_{s}\neq0\}}
|Z^{n}_{s}-Z^{m}_{s}|^{2}\,ds\\
&\nonumber\quad\le c(p,\lambda) E\int_{0}^{T}|Y^{n}_{s}-Y^{m}_{s}|^{p}\,ds
+(E|(Y^{n}-L)^{-,*}_{T}|^{p})^{(p-1)/p}(E|K^{m}_{T}|^{p})^{1/p}
\\&\qquad
+(E|(Y^{m}-L)^{-,*}_{T}|^{p})^{(p-1)/p}(E|K^{n}_{T}|^{p})^{1/p}\equiv I_{n,m}.
\end{align*}
From the above, (\ref{eq4.2}), (\ref{eq4.4}) and the monotone
convergence of $\{Y^{n}\}$ we get
\begin{equation}
\label{T3}
\lim_{n,m\rightarrow +\infty} I^{n,m}=0
\end{equation}
which implies that
\begin{equation}\label{eq4.9}
\lim_{n,m\rightarrow+\infty}E\int_{0}^{T}|Y^{n}_{s}-Y^{m}_{s}|^{p-2}
\mathbf{1}_{\{Y^{n}_{s}\neq Y^{m}_{s}\}}
|Z^{n}_{s}-Z^{m}_{s}|^{2}\,ds=0.
\end{equation}
From (\ref{eq4.5}) one can also conclude that
\begin{align*}
&E\sup_{0\le t\le T}|Y^{n}_{t}-Y^{m}_{t}|^{p}\\
&\qquad\le c'(p,\lambda) \{I^{n,m}+E\sup_{0\le t\le T}|\intt
|Y^{n}_{s}-Y^{m}_{s}|^{p-1}
\widehat{Y^{n}_{s}-Y^{m}_{s}}(Z^{n}_{s}-Z^{m}_{s})\,dB_{s}|\}.
\end{align*}
Using the Burkholder-Davis-Gundy inequality and then Young's
inequality we deduce from the above that
\[
E\sup_{0\le t\le T}|Y^{n}_{t}-Y^{m}_{t}|^{p}\le c''(p,\lambda)
\{I^{n,m}+E\int_{0}^{T}\mathbf{1}_{\{Y^{n}_{s}
\neq Y^{m}_{s}\}} |Y^{n}_{s}-Y^{m}_{s}|^{p-2}
|Z^{n}_{s}-Z^{m}_{s}|^{2}\,ds\}.
\]
Hence, by (\ref{T3}) and (\ref{eq4.9}),
\begin{align}\label{eq4.10}
\lim_{n,m\rightarrow+\infty}E\sup_{0\le t\le
T}|Y^{n}_{t}-Y^{m}_{t}|^{p}=0,
\end{align}
which implies that $Y\in\mathcal{S}^{p}$. Our next goal is to show
that
\begin{equation}\label{eq4.11}
\lim_{n,m\rightarrow+\infty}
E(\int_{0}^{T}|Z^{n}_{t}-Z^{m}_{t}|^{2}\,dt)^{p/2}=0.
\end{equation}
By It\^o's formula applied to $|Y^{n}-Y^{m}|^{2}$, (H2) and (H3),
\begin{align*}
\int_{0}^{T}|Z^{n}_{t}-Z^{m}_{t}|^{2}\,dt &\le
2\lambda\int_{0}^{T}|Y^{n}_{t}-Y^{m}_{t}||Z^{n}_{t}-Z^{m}_{t}|\,dt
+2\int_{0}^{T}|Y^{n}_{t}-Y^{m}_{t}|\,dK^{n}_{t}\\&\quad
+2\int_{0}^{T}|Y^{n}_{t}-Y^{m}_{t}|\,dK^{m}_{t} +\sup_{0\le t\le
T}|\intt (Z^{n}_{s}-Z^{m}_{s}) (Y^{n}_{s}-Y^{m}_{s})\,dB_{s}|.
\end{align*}
Hence, by the Burkholder-Davis-Gundy inequality and  Young's
inequality,
\begin{align*}
&E(\int_{0}^{T}|Z^{n}_{t}-Z^{m}_{t}|^{2}\,dt)^{p/2}\le
C(p,\lambda)\{E|(Y^{n}-Y^{m})^{*}_{T}|^{p}\\
&\qquad+(E|(Y^{n}-Y^{m})^{*}_{T}|^{p})^{1/2}(E|K^{n}_{T}|^{p})^{1/2}
+(E|(Y^{n}-Y^{m})^{*}_{T}|^{p})^{1/2}(E|K^{m}_{T}|^{p})^{1/2}\}.
\end{align*}
From the above inequality, (\ref{eq4.2}) and (\ref{eq4.10}) we get
(\ref{eq4.11}). From (\ref{eq4.11}) and (\ref{eq4.3}) it follows
immediately that
\[
Y_{t}=\xi+\int_{t}^{T}f(s,Y_{s},Z_{s})\,ds+\int_{t}^{T}dK_{s}
-\int_{t}^{T}Z_{s}\,dB_{s},\quad t\in [0,T],
\]
which implies that $K$ is continuous. In fact, by (\ref{eq4.2}),
$K\in \mathcal{V}^{+,p}_{c}$. Moreover, from (\ref{eq5.01}),
(\ref{T1}), (\ref{T2}) (\ref{eq4.10}), (\ref{eq4.11}) and (H2) we
deduce that
\begin{align}\label{eq4.12}
\lim_{n,m\rightarrow+\infty}E\sup_{0\le t\le T}
|K^{n}_{t}-K^{m}_{t}|^{p}=0.
\end{align}
Since $\int_{0}^{T}(Y^{n}_{t}-L_{t})\,dK^{n}_{t}\le 0$, it follows
from  (\ref{eq4.10}), (\ref{eq4.12}) that
$\int_{0}^{T}(Y_{t}-L_{t})\,dK_{t}\le0$, which when combined with
the fact that $Y_{t}\ge L_{t}$, $t\in [0,T]$ shows that
\[
\int_{0}^{T}(Y_{t}-L_{t})\,dK_{t}=0.
\]
Thus the triple $(Y,Z,K)$ is a solution of RBSDE$(\xi,f,L)$, which
completes the proof of (i). Assertion (ii) follows from
(\ref{eq4.10})--(\ref{eq4.12}).
\end{dow}

\begin{uw}
Let $p>1$ and let assumptions (H1)--(H3) hold. If $(Y,Z,K)$ is a
solution of RBSDE$(\xi,f,L)$ such that $(Y,Z)\in
\mathcal{S}^{p}\otimes M^{p}$ then from \cite[Remark 4.3]{BDHPS}
it follows immediately that
\[
E(\int_{0}^{T}|f(s,Y_{s},Z_{s})|\,ds)^{p}<+\infty \mbox{ iff }
EK^{p}_{T}<+\infty.
\]
Moreover, if there exists $X\in\HH^{p}_{c}$ such that
$E(\int_{0}^{T}f^{-}(s,X_{s},0)\,ds)^{p}<+\infty$ then
\begin{equation}
\label{eq5.14} E(\int_{0}^{T}\mathbf{1}_{\{Y_{s}\le X_{s}\}}\,
dK_{s})^{p}<+\infty.
\end{equation}
Indeed, since $X\in\HH^{p}_{c} $, there exist $M\in\MM^{p}_{c}$
and $V\in\mathcal{V}^{p}_{c}$ such that $X_{t}=X_{0}+M_{t}+V_{t}$,
$t\in [0,T]$. Let $L^{0}(Y-X)$ denote the local time of $Y-X$ at
0. By (H2), (H3) and the It\^o-Tanaka formula applied to
$(Y-X)^{-}$,
\begin{align*}
\int_{0}^{T}\mathbf{1}_{\{Y_{s}\le X_{s}\}}\,dK_{s}
&=(Y_{T}-X_{T})^{-}-(Y_{0}-X_{0})^{-}
-\int_{0}^{T}\mathbf{1}_{\{Y_{s}\le X_{s}\}}f(s,Y_{s},Z_{s})\,ds
\\&\quad
-\int_{0}^{T}\mathbf{1}_{\{Y_{s}\le X_{s}\}} dV_{s}
-\frac12\int_{0}^{T}dL^{0}_{s}(Y-X)
-\int_{0}^{T}\mathbf{1}_{\{Y_{s}\le X_{s}\}} Z_{s}\,dB_{s}\\
&\quad +\int_{0}^{T}\mathbf{1}_{\{Y_{s}\le X_{s}\}}\,dM_{s}\\
&\le2Y^{*}_{T}+2X^{*}_{T}-\int_{0}^{T}\mathbf{1}_{\{Y_{s} \le
X_{s}\}}f(s,X_{s},0)\,ds+\lambda\int_{0}^{T}|Z_{s}|\,ds\\
&\quad+\int_{0}^{T}d|V|_{s}-\int_{0}^{T}\mathbf{1}_{\{Y_{s}\le
X_{s}\}} Z_{s}\,dB_{s}+\int_{0}^{T}\mathbf{1}_{\{Y_{s}\le X_{s}\}}
\,dM_{s},
\end{align*}
from which one can easily get (\ref{eq5.14}).
\end{uw}

We close this section with an example  which shows that assumption
(\ref{i7}) is not necessary for existence of $p$-integrable
solutions of reflected BSDEs.

\begin{prz}
Let $V_{t}=\exp(|B_{t}|^{4})$, $t\in [0,T]$. Observe that
\[
P(\int_{0}^{T}V_{t}\,dt<+\infty)=1,\quad E\int_{a}^{T}V_{t}\,dt
=+\infty,\quad a\in(0,T).
\]
Now, set $\xi\equiv 0$, $f(t,y)=-(y-(T-t))^{+}V_{t}$, $L_{t}=T-t$,
$t\in [0,T]$. Then  $\xi,f,L$ satisfy (H1)--(H7) with $p=2$. On
the other hand,
\begin{align*}
E\int_{0}^{T}f^{-}(t,L^{*}_{t})\,dt=E\int_{0}^{T}f^{-}(t,T)\,dt
=E\int_{0}^{T}tV_{t}\,dt\ge aE\int_{a}^{T}V_{t}\,dt=+\infty.
\end{align*}
\end{prz}

\nsubsection{Existence and uniqueness results for $p=1$}
\label{sec6}

We first prove uniqueness.
\begin{stw}\label{stw4}
If $f$ satisfies \mbox{\rm(H2)}, \mbox{\rm(H3)} and \mbox{\rm(Z)}
then there exists at most one solution $(Y,Z,K)$ of \mbox{\rm
RBSDE}$(\xi,f,L)$ such that $Y$ is of class \mbox{\rm(D)} and
$Z\in\bigcup_{\beta>\alpha} M^{\beta}$.
\end{stw}
\begin{dow}
Without loss of generality we may assume that $\mu\le 0$. Let
$(Y^{1},Z^{1},K^{1})$, $(Y^{2},Z^{2},K^{2})$ be two solutions  to
RBSDE$(\xi,f,L)$. By Proposition \ref{stw2.5} it suffices to prove
that $|Y^{1}-Y^{2}|\in \mathcal{S}^{p}$ for some $p>1$. Write
$Y=Y^{1}-Y^{2}$, $Z=Z^{1}-Z^{2}$, $K=K^{1}-K^{2}$ and
$\tau_{k}=\inf\{t\in[0,T];
\int_{0}^{t}(|Z^{1}_{s}|^{2}+|Z^{2}_{s}|^{2})\,ds>k\}\wedge T$.
Then by the It\^o formula (see \cite[Corollary 2.3]{BDHPS}),
\begin{align*}
|Y_{t\wedge\tau_{k}}|&\le |Y_{\tau_{k}}|
+\int_{t\wedge \tau_{k}}^{\tau_{k}}\hat{Y}_{s}
(f(s,Y^{1}_{s},Z^{1}_{s})-f(s,Y^{2}_{s},Z^{2}_{s}))\,ds\\&\quad
+\int_{t\wedge \tau_{k}}^{\tau_{k}}\hat{Y}_{s}\,dK_{s}
-\int_{t\wedge \tau_{k}}^{\tau_{k}}\hat{Y}_{s}Z_{s}\,dB_{s},
\quad t\in [0,T].
\end{align*}
By the minimality property (d) of the reaction measures $K^{1},
K^{2}$ in the definition of a solution of RBSDE$(\xi,f,L)$,
$\int_{0}^{T}\hat{Y}_{s}\,dK_{s}\le 0$. Hence
\begin{align*}
|Y_{t\wedge\tau_{k}}|&\le |Y_{\tau_{k}}| +\int_{t\wedge
\tau_{k}}^{\tau_{k}}\hat{Y}_{s}
(f(s,Y^{1}_{s},Z^{1}_{s})-f(s,Y^{2}_{s},Z^{2}_{s}))\,ds
-\int_{t\wedge \tau_{k}}^{\tau_{k}} \hat{Y}_{s}Z_{s}\,dB_{s}\\
&\le |Y_{\tau_{k}}|
+\int_{0}^{T}|f(s,Y^{1}_{s},Z^{1}_{s})-f(s,Y^{1}_{s},Z^{2}_{s})|\,ds
-\int_{t\wedge \tau_{k}}^{\tau_{k}}\hat{Y}_{s}Z_{s}\,dB_{s}
\end{align*}
for $t\in[0,T]$, the last inequality being a consequence of (H3).
Consequently,
\begin{align*}
|Y_{t\wedge\tau_{k}}|\le E^{\FF_{t}}(|Y_{\tau_{k}}|
+\int_{0}^{T}|f(s,Y^{1}_{s},Z^{1}_{s})-f(s,Y^{1}_{s},Z^{2}_{s})|\,ds),
\quad t\in [0,T].
\end{align*}
Since $Y$ is of class (D), letting $k\rightarrow +\infty$ we conclude from
the above that
\begin{align*}
|Y_{t}|\le E^{\FF_{t}}(\int_{0}^{T}
|f(s,Y^{1}_{s},Z^{1}_{s})-f(s,Y^{1}_{s},Z^{2}_{s})|\,ds),\quad t\in [0,T].
\end{align*}
By (Z),
\[
|Y_{t}|\le 2\gamma E^{\FF_{t}}(\int_{0}^{T}
(g_{s}+|Y^{1}|_{s}+|Z^{1}|_{s}+|Z^{2}_{s}|)^{\alpha}\,ds).
\]
From this it follows that  $|Y|\in \mathcal{S}^{p}$ for some
$p>1$, which proves the proposition.
\end{dow}

\begin{uw}\label{uw.dic}
A brief inspection of the proof of Proposition \ref{stw4} reveals
that if $f$ does not depend on $z$ and satisfies (H2) then there
exits at most one solution $(Y,Z,K)$ of RBSDE$(\xi,f,L)$ such that
$Y$  is of class (D).
\end{uw}

\begin{uw}\label{uw2}
If (H1), (H3), (Z) are satisfied and $(Y,Z)$ is a unique solution
of BSDE$(\xi,f)$ such that $Y$ is of class (D) and
$Z\in\mathbb{L}^{\alpha}(\FF)$ then
\[
E\int_{0}^{T}|f(s,Y_{s},Z_{s})|\,ds<+\infty.
\]
Indeed, by Proposition \ref{prop.ito}, for every stopping time
$\tau\le T$,
\begin{align*}
|Y_{t\wedge\tau}|\le |Y_{\tau}|
+\int_{t\wedge\tau}^{\tau}\hat{Y}_{s}f(s,Y_{s},Z_{s})\,ds
-\int_{t\wedge\tau}^{\tau}\hat{Y}_{s}Z_{s}\,dB_{s},\quad t\in [0,T].
\end{align*}
Hence
\begin{align*}
-\int_{t\wedge\tau}^{\tau}\hat{Y}_{s}
(f(s,Y_{s},Z_{s})-f(s,0,Z_{s}))\,ds &\le
|Y_{\tau}|-|Y_{t\wedge\tau}|
+\int_{t\wedge\tau}^{\tau}|f(s,0,Z_{s})|\,ds\\&\quad
-\int_{t\wedge\tau}^{\tau}\hat{Y}_{s}Z_{s}\,dB_{s}.
\end{align*}
By the above inequality, (H3) (without loss of generality we may
assume that $\mu\le 0$) and (Z), for $t\in[0,T]$ we have
\begin{align*}
&E\int_{t\wedge\tau_k}^{\tau}|f(s,Y_{s},Z_{s})-f(s,0,Z_{s})|\,ds\\
&\qquad\le E|Y_{\tau_k}|
+E\int_{t\wedge\tau_k}^{\tau}(g_{s}+|Z_{s}|+|Y_{s}|)^{\alpha}\,ds
+\int_{t\wedge\tau_k}^{\tau}f_{s}\,ds,
\end{align*}
where $\tau_{k}$ is defined by (\ref{eq3.01}). Since $Y$  is of class (D),
letting $k\rightarrow +\infty$ we obtain
\begin{align*}
E\int_{0}^{T}|f(s,Y_{s},Z_{s})-f(s,0,Z_{s})|\,ds \le E|\xi|+\gamma
E\int_{0}^{T}(g_{s}+|Z_{s}|+|Y_{s}|)^{\alpha}\,ds
+\int_{0}^{T}f_{s}\,ds.
\end{align*}
Using once again (Z) we conclude from the above that
\begin{align*}
E\int_{0}^{T}|f(s,Y_{s},Z_{s})|\,ds \le E|\xi|+2\gamma
E\int_{0}^{T}(g_{s}+|Z_{s}|+|Y_{s}|)^{\alpha}\,ds+ 2\int_{0}^{T}
f_{s}\,ds<+\infty.
\end{align*}
\end{uw}

\begin{tw}\label{tw2}
Let $p=1$.
\begin{enumerate}
\item[\rm(i)]Assume \mbox{\rm(H1)}--\mbox{\rm(H6)}, \mbox{\rm(Z)}. Then
there exists a solution $(Y,Z,K)$ of \mbox{\rm RBSDE}$(\xi,f,L)$
such that $Y$  is of class \mbox{\rm(D)}, $K\in
\mathcal{V}^{+,1}_{c}$ and $Z\in\bigcap_{q<1}M^{q}$ iff
\mbox{\rm(H7*)} is satisfied.
\item[\rm(ii)]Assume \mbox{\rm(H1)}--\mbox{\rm(H6)}, \mbox{\rm(H7*)} and
for $n\in\BN$ let $(Y^{n},Z^{n})$ be a solution of
\mbox{\rm(\ref{eq5.01})} such that $(Y^{n},Z^{n})\in
\mathcal{S}^{q}\otimes M^{q}$, $q\in(0,1)$, and $Y^{n}$ is of
class \mbox{\rm(D)}. Let  $K^{n}$ be  defined by
\mbox{\rm(\ref{eq5.02})}. Then for every $q\in(0,1)$,
\[
E\sup_{t\le T}|Y^{n}_{t}-Y_{t}|^{q} +E\sup_{t\le
T}|K^{n}_{t}-K_{t}|^{q}
+E(\int_{0}^{T}|Z^{n}_{t}-Z_{t}|^{2}\,dt)^{q/2}\rightarrow 0
\]
as $n\rightarrow +\infty$.
\end{enumerate}
\end{tw}
\begin{dow}
(i) Necessity. By Remark \ref{uw2}, if there is a solution
$(Y,Z,K)$ of BSDE$(\xi,f,L)$ such that $(Y,Z)\in
\mathcal{S}^{q}\otimes M^{q}$, $q\in (0,1)$,
$K\in\mathcal{V}^{+,1}_{c}$ and $Y$ is of class (D) then (H7*) is
satisfied with $X=Y$.
\smallskip\\
Sufficiency. We first show that the sequence $\{Y^{n}\}$ is
nondecreasing. To this end, let us put
$f_{n}(t,y,z)=f(t,y,z)+n(y-L_{t})^{-}$. Since the exponential
change of variable described at the beginning of the proof of
Lemma \ref{lm1} does not change the monotonicity of the sequence
$\{Y^{n}\}$, we may and will assume that the mapping
$\mathbb{R}\ni y\mapsto f_{n}(t,y,0)$ is nonincreasing. By the
It\^o-Tanaka formula, for every stopping time $\tau\le T$,
\begin{align*}
&(Y^{n}_{t\wedge\tau}-Y^{n+1}_{t\wedge\tau})^{+}
+\frac12\int_{\tau\wedge t}^{\tau}dL^{0}_{s}(Y^{n}-Y^{n+1})\\
&\qquad=(Y^{n}_{\tau}-Y^{n+1}_{\tau})^{+}
+\int_{t\wedge\tau}^{\tau} \mathbf{1}_{\{Y^{n}_{s} >
Y^{n+1}_{s}\}}(f_{n}(s,Y^{n}_{s},Z^{n}_{s})
-f_{n+1}(s,Y^{n+1}_{s},Z^{n+1}_{s}))\,ds\\
&\qquad\quad-\int_{t\wedge\tau}^{\tau} \mathbf{1}_{\{Y^{n}_{s} >
Y^{n+1}_{s}\}}(Z^{n}_{s}-Z^{n+1}_{s})\,dB_{s}.
\end{align*}
Taking the conditional expectation with respect to $\FF_t$ on both
sides of the above equality with $\tau$ replaced by
$\tau_k=\inf\{t\in [0,T];\int_{0}^{t}
|Z^{n}_{s}-Z^{n+1}_{s}|^{2}\,ds\ge k\}\wedge T$,  letting
$k\rightarrow +\infty$ and using the fact that $Y$  is of class (D) we obtain
\begin{align}
\label{T4} (Y^{n}_{t}-Y^{n+1}_{t})^{+}&\le
E^{\mathcal{F}_{t}}\int_{t}^{T}
\mathbf{1}_{\{Y^{n}_{s}>Y^{n+1}_{s}\}}(f_{n}(s,Y^{n}_{s},Z^{n}_{s})
-f_{n+1}(s,Y^{n+1}_{s},Z^{n+1}_{s}))\,ds.
\end{align}
From the above inequality and the fact that $f_{n}\le
f_{n+1}$ we get
\begin{align*}
&\int_{t}^{T} \mathbf{1}_{\{Y^{n}_{s}
> Y^{n+1}_{s}\}}(f_{n}(s,Y^{n}_{s},Z^{n}_{s})
-f_{n+1}(s,Y^{n}_{s},Z^{n}_{s}))\,ds\\&\quad
+\int_{t}^{T} \mathbf{1}_{\{Y^{n}_{s}
> Y^{n+1}_{s}\}}(f_{n+1}(s,Y^{n}_{s},Z^{n}_{s})
-f_{n+1}(s,Y^{n+1}_{s},Z^{n+1}_{s}))\,ds
\\&\le \int_{t}^{T} \mathbf{1}_{\{Y^{n}_{s}
>Y^{n+1}_{s}\}}(f_{n+1}(s,Y^{n}_{s},Z^{n}_{s})
-f_{n+1}(s,Y^{n+1}_{s},Z^{n+1}_{s}))\,ds\\&=
\int_{t}^{T}\mathbf{1}_{\{Y^{n}_{s} >
Y^{n+1}_{s}\}}(f_{n+1}(s,Y^{n}_{s},Z^{n}_{s})
-f_{n+1}(s,Y^{n}_{s},0))\,ds
\\&\quad+\int_{t}^{T}\mathbf{1}_{\{Y^{n}_{s}
> Y^{n+1}_{s}\}}(f_{n+1}(s,Y^{n}_{s},0)
-f_{n+1}(s,Y^{n+1}_{s},0))\,ds
\\&\quad+\int_{t}^{T}\mathbf{1}_{\{Y^{n}_{s}
> Y^{n+1}_{s}\}}(f_{n+1}(s,Y^{n+1}_{s},0)
-f_{n+1}(s,Y^{n+1}_{s},Z^{n+1}_{s}))\,ds.
\end{align*}
Since $f_{n}(t,y,z)-f_{n}(t,y,z')=f(t,y,z)-f(t,y,z')$ for every
$t\in[0,T]$, $y\in\mathbb{R}$, $z,z'\in\mathbb{R}^{d}$, using the
monotonicity of $f_{n+1}$ and assumption (Z) we conclude from the
above and (\ref{T4}) that for $t\in[0,T]$,
\[
(Y^{n}_{t}-Y^{n+1}_{t})^{+} \le 2\gamma
E^{\mathcal{F}_{t}}\int_{0}^{T}
(g_{s}+|Y^{n}_{s}|+|Z_{s}^{n}|+|Y^{n+1}_s|+|Z^{n+1}_s|)^{\alpha}\,ds.
\]
Since $(Y^{n},Z^{n})\in\mathcal{S}^{q}\otimes M^{q}$ for every
$q\in (0,1)$, $n\in\mathbb{N}$,  it follows from the above
estimate that $(Y^{n}-Y^{n+1})^{+}\in \mathcal{S}^{p}$ for some
$p>1$. Hence, by Proposition \ref{stw2}, $Y^{n}_{t}\le
Y^{n+1}_{t}$, $t\in [0,T]$. Write
\[
Y_{t}=\lim_{n\rightarrow +\infty} Y^{n}_{t},\quad t\in [0,T].
\]
We are going to show that there is a process $\overline{X}$ of
class (D) such that $\overline{X}\in
\mathcal{V}^{1}_{c}+\MM^{q}_{c}$ for $q\in(0,1)$ and
$\overline{X}_{t}\ge Y_{t}$, $t\in [0,T]$. Indeed, since $X$ from
assumption (H7*) belongs to $\mathcal{V}^{1}_{c}+\MM^{q}_{c}$ for
$q\in (0,1)$, there exist $M\in\MM^{q}_{c}$ and
$V\in\mathcal{V}^{1}_{c}$ such that $X=V+M$. By the representation
property of the Brownian filtration there exists $Z'\in M^{q}$
such that
\[
X_{t}=X_{T}-\intt dV_{s}-\intt Z'_{s}\,dB_{s},\quad t\in [0,T],
\]
which we can write in the form
\begin{align*}
X_{t}&=X_{T}+\intt f(s,X_{s},Z'_{s})\,ds
-\intt (f^{+}(s,X_{s},Z'_{s})\,ds+dV^{+}_{s})\\&\quad
+\intt (f^{-}(s,X_{s},Z'_{s})\,ds+dV^{-}_{s})
-\intt Z'_{s}\,dB_{s},\quad t\in [0,T].
\end{align*}
By \cite[Theorem 6.3]{BDHPS} and Remark \ref{uw2} there exists a
unique solution $(\overline{X},\overline{Z})$ of the BSDE
\begin{align*}
\overline{X}_{t}=\xi\vee X_T+\intt
f(s,\overline{X}_{s},\overline{Z}_{s})\,ds +\intt
(f^{-}(s,X_{s},Z'_{s})\,ds+dV^{-}_{s}) -\intt
\overline{Z}_{s}\,dB_{s}
\end{align*}
such that
$(\overline{X},\overline{Z})\in\bigcap_{q<1}\mathcal{S}^{q}
\otimes M^{q}$, $\overline{X}$ is of class (D) and
\begin{equation}
\label{T5}
E\int_{0}^{T}|f(t,\bar{X}_{t},\bar{Z}_{t})|\,dt<+\infty.
\end{equation}
As in the proof of the fact that $(Y^{n}-Y^{n+1})^{+}\in
\mathcal{S}^{p}$ one can show that for every stopping time
$\tau\le T$,
\begin{align*}
(X_{t\wedge\tau}-\overline{X}_{t\wedge\tau})^{+}
&\le(X_{\tau}-\overline{X}_{\tau})^{+}
+\int_{t\wedge\tau}^{\tau}\mathbf{1}_{\{X_{s}>\overline{X}_{s}\}}
(f(s,X_{s},Z_{s}')-f(s,\overline{X}_{s},\overline{Z}_{s}))\,ds\\
&\quad-2\int_{t\wedge\tau}^{\tau}\mathbf{1}_{\{X_{s}>\overline{X}_{s}\}}
(Z'_{s}-\overline{Z}_{s})\,dB_{s}\\
&\le (X_{\tau}-\overline{X}_{\tau})^{+}
+2\gamma\int_{t\wedge\tau}^{\tau}(g_{s}+|X_{s}|+|\bar{X}_{s}|+|Z'_{s}|+|Z_{s}|)^{\alpha}\,ds\\
&\quad-2\int_{t\wedge\tau}^{\tau}\mathbf{1}_{\{X_{s}
>\overline{X}_{s}\}}(Z'_{s}-\overline{Z}_{s})\,dB_{s}.
\end{align*}
Let $\tau_{k}=\inf\{t\in [0,T];\int_{0}^{t}
(|Z^{'}_{s}|^{2}+|\overline{Z}_{s}|^{2})\,ds\ge k\}\wedge T$. Then
\[
(X_{t\wedge\tau_{k}}-\overline{X}_{t\wedge\tau_{k}})^{+}
\le E^{\mathcal{F}_{t}}(X_{\tau_{k}}-\overline{X}_{\tau_{k}})^{+}
+2\gamma E^{\FF_{t}}\int_{0}^{T}(g_{s}+|X_{s}|+|\bar{X}_{s}|+|Z'_{s}|+|Z_{s}|)^{\alpha}\,ds.
\]
Since $X,\overline{X}$ are of class (D), letting
$k\rightarrow+\infty$ we get
\[
(X_{t}-\overline{X}_{t})^{+}\le 2\gamma
E^{\FF_{t}}\int_{0}^{T}
(g_{s}+|X_{s}|+|\bar{X}_{s}|+|Z'_{s}|+|Z_{s}|)^{\alpha}\,ds.
\]
Therefore $(X-\overline{X})^{+}\in \mathcal{S}^{p}$ for some $p>1$
since $Z',\overline{Z}\in M^{q}$, $X,\bar{X}\in\mathcal{S}^{q}$,
$q\in (0,1)$. Consequently, by Proposition \ref{stw2}, $X_{t}\le
\overline{X}_{t}$, $t\in [0,T]$. Thus,
\begin{align*}
\overline{X}_{t}&=\xi\vee X_{T}
+\intt f(s,\overline{X}_{s},\overline{Z}_{s})\,ds
+\intt n(\overline{X}_{s}-L_{s})^{-}\,ds
\\&\quad+\intt (f^{-}(s,X_{s},Z'_{s})\,ds
+dV^{-}_{s})-\intt \overline{Z}_{s}\,dB_{s},\quad t\in [0,T].
\end{align*}
As in the case of the process $(X-\overline{X})^{+}$ one can show
that $(Y^{n}-\overline{X})^{+}\in\mathcal{S}^{p}$ for some $p>1$.
Hence, by Proposition \ref{stw2},  $Y^{n}_{t}\le
\overline{X}_{t}$, $t\in [0,T]$ for every $n\in\mathbb{N}$.
Furthermore, since $Y^{1},\overline{X}\in \mathcal{S}^{q}$, $q\in
(0,1)$, we have
\begin{equation}
\label{eq5.1} \sup_{n\ge1}E|Y^{n,*}_{T}|^{q}<+\infty.
\end{equation}
It follows in particular that $\sup_{n\ge1}|Y^{n}_{0}|<\infty$
since $Y^{n}_{0}$ are deterministic. Moreover, by Lemma \ref{lm4},
there exists a stationary sequence $\{\sigma^{1}_{k}\}$ of
stopping times such that for every $k\in\mathbb{N}$,
\begin{equation}
\label{eq5.2} \sup_{n\ge1}|Y^{n,*}_{\sigma^{1}_{k}}| \le
k\vee(\sup_{n\ge1}|Y^{n}_{0}|)<+\infty.
\end{equation}
Set
\[
\sigma^{2}_{k}=\inf\{t\in [0,T],
\min\{Y^{1,*}_{t},\overline{X}^{+,*}_{t},
\int_{0}^{t}f^{-}(s,\overline{X}_{s},0)\,ds,
\int_{0}^{t}|f(s,0,0)|\,ds\}>k\}\wedge T
\]
and $\tau_{k}=\sigma^{1}_{k}\wedge\sigma^{2}_{k}$. It is easy to
see that the sequence $\{\tau_{k}\}$ is stationary. Using this and
the fact that $Y^n_{\tau_k}$, $f$, $L$ satisfy the assumptions of
Theorem \ref{tw1} on the interval $[0,\tau_k]$ one can show that there
exist $Y,K\in\mathcal{S}, Z\in M$ such that $K$ is increasing,
$K_{0}=0$ and
\begin{equation}
\label{eq5.3}
\sup_{0\le t \le T}|Y^{n}_{t}-Y_{t}| +\sup_{0\le t\le T}
|K^{n}_{t}-K_{t}|+\int_{0}^{T}|Z^{n}_{s}-Z_{s}|^{2}\,ds\rightarrow
0\mbox{ in probability }P
\end{equation}
as $n\rightarrow +\infty$. Moreover, one can show that  $Y_{t}\ge
L_{t}$, $t\in [0,T]$,
\begin{align}\label{eq5.4}
Y_{t}=\xi+\int_{t}^{T}f(s,Y_{s},Z_{s})\,ds
+\int_{t}^{T}dK_{s}-\int_{t}^{T} Z_{s}\,dB_{s},\quad t\in[0,T]
\end{align}
and
\begin{equation}
\label{eq5.5} \int_{0}^{T}(Y_{s}-L_{s})\,dK_{s}=0.
\end{equation}
Accordingly, the trip $(Y,Z,K)$ is a solution of RBSDE$(\xi,f,L)$.
The proof of (\ref{eq5.3})--(\ref{eq5.5}) runs as the proof of
Theorem \ref{tw1} (see the reasoning following (\ref{eq4.2}) with
$p=2$), the only difference being in the fact that now we consider
equations on $[0,\tau_{k}]$ with terminal values depending on $n$.
However, using (\ref{eq5.2}) and the pointwise convergence of
$\{Y^{n}\}$ allows overcome this difficulty. Since $Y^{1}_{t}\le
Y_{t}\le \overline{X}_{t}$, $t\in [0,T]$, and $Y^{1}, X^{+}$  are
of class (D), it follows that $Y$  is of class (D). By Lemma
\ref{lm2} for every $q\in(0,1)$,
\begin{equation}
\label{eq6.8}
\sup_{n\ge1}E\big((\int_{0}^{T}|Z^{n}_{t}|^{2}\,dt)^{q/2}
+|K^{n}_{T}|^{q}\big)<+\infty.
\end{equation}
From this and (\ref{eq5.3}) we conclude that
$Z\in\bigcap_{q<1}M^{q}$ and $E|K_{T}|^{q}<\infty$ for
$q\in(0,1)$. To see that $EK_{T}<\infty$ let us define $\tau_{k}$
by (\ref{eq3.01}). Then by (\ref{eq5.4}),
\begin{equation}
\label{eq6.9} K_{\tau_{k}}=Y_{0}-Y_{\tau_{k}}
-\int_{0}^{\tau_{k}}f(s,Y_{s},Z_{s})\,ds +\int_{0}^{\tau_{k}}
Z_{s}\,dB_{s}.
\end{equation}
Since $Y$  is of class (D), using Fatou's lemma, (H2), (Z) and the
fact that $Y_{t}\le\overline{X}_{t}$, $t\in [0,T]$ we conclude
from (\ref{eq6.9}) that
\begin{align*}
EK_{T}\le EY^{+}_{0}+E\xi^{-}
+E\int_{0}^{T}f^{-}(s,\overline{X}_{s},0)\,ds+\gamma
E\int_{0}^{T}(g_{s}+|Y_{s}|+|Z_{s}|)^{\alpha}\,ds.
\end{align*}
Hence $EK_{T}<\infty$, because by (\ref{T5}) and (H2),
$E\int_{0}^{T}|f(s,\overline{X}_{s},0)|\,ds<+\infty$.
\smallskip\\
(ii) Convergence of $\{Y^{n}\}$ in $\mathcal{S}^{q}$ for $q\in
(0,1)$ follows from (\ref{eq5.1}) and (\ref{eq5.3}). The desired
convergence of $\{Z^{n}\}$ and $\{K^{n}\}$ follows from
(\ref{eq5.3}) and (\ref{eq6.8}).
\end{dow}

\begin{uw}
An important class of generators satisfying (H1)--(H5) together
with (Z) are generators satisfying (H1)--(H5) which are bounded or not
depending on $z$. Another class which share these properties are
generators of the form
\[
f(t,y,z)=g(t,y)+c(1+|z|)^{q},
\]
where $q\in [0,\alpha]$ and $g$ is a progressively measurable
function satisfying (H1)--(H5).
\end{uw}

\begin{uw}
Let assumptions (H1)--(H3), (Z) hold and let $(Y,Z,K)$ be a
solution of RBSDE$(\xi,f,L)$ such that $Y$  is of class (D) and
$Z\in\bigcup_{\beta>\alpha} M^{\beta}$. Then from Remark \ref{uw1}
it follows immediately that
\[
E(\int_{0}^{T}|f(s,Y_{s},Z_{s})|\,ds)<+\infty\mbox{ iff }
EK_{T}<+\infty.
\]
If, in addition, there exists a continuous semimartingale $X$ such
that (H7*) is satisfied then
\[
E\int_{0}^{T}\mathbf{1}_{\{Y_{s}\le X_{s}\}}dK_{s}<+\infty.
\]
To prove the last estimate let us put $\tau_{k}=\inf\{t\in [0,T];
\langle M\rangle_{t} +\int_{0}^{t}|Z_{s}|^{2}\,ds>k\}\wedge T$. By
the It\^o-Tanaka formula and (H2), (H3),
\begin{align*}
\int_{0}^{\tau_{k}}\mathbf{1}_{\{Y_{s}\le X_{s}\}}\,dK_{s}
&=(Y_{\tau_{k}}-X_{\tau_{k}})^{-}-(Y_{0}-X_{0})^{-}
-\int_{0}^{\tau_{k}}\mathbf{1}_{\{Y_{s}\le X_{s}\}}
f(s,Y_{s},Z_{s})\,ds\\
&\quad-\int_{0}^{\tau_{k}}\mathbf{1}_{\{Y_{s}\le X_{s}\}}dV_{s}
-\frac12\int_{0}^{\tau_{k}}dL^{0}_{s}(Y-X)\\
&\quad-\int_{0}^{\tau_{k}}\mathbf{1}_{\{Y_{s}\le X_{s}\}}
Z_{s}\,dB_{s}+\int_{0}^{\tau_{k}}\mathbf{1}_{\{Y_{s}\le
X_{s}\}}\,dM_{s}.
\end{align*}
Hence
\begin{align*}
E\int_{0}^{\tau_{k}}\mathbf{1}_{\{Y_{s}\le X_{s}\}}\,dK_{s} &\le
E|Y_{\tau_{k}}|+ EX^{+}_{\tau_{k}}
+E\int_{0}^{T}\mathbf{1}_{\{Y_{s}\le X_{s}\}}
f^{-}(s,X_{s},0)\,ds\\&\quad +\gamma
E\int_{0}^{T}(g_{s}+|Z_{s}|+|Y_{s}|)^{\alpha}\,ds
+E\int_{0}^{T}d|V|_{s}.
\end{align*}
Since $(Y-X)^{-}$ is of class (D), letting $k\rightarrow+\infty$
in the above inequality we get the desired result.
\end{uw}

\nsubsection{Nonintegrable solutions of reflected BSDEs}
\label{sec7}

In this section we examine existence and uniqueness of solutions
of reflected BSDEs in the case where the data satisfy (H1)--(H6)
(resp. (H1)--(H6), (Z) for $p=1$)  but (H7) (resp. (H7*) in case
$p=1$) is not satisfied. In view of Theorems \ref{tw1} and
\ref{tw2} in that case there is neither a solution $(Y,Z,K)$ in
the space $\mathcal{S}^{p}\otimes M^{p}\otimes
\mathcal{V}^{p,+}_{c}$ if $p>1$ nor a solution in the space
$\mathcal{S}^{q}\otimes M^{q}\otimes \mathcal{V}^{1,+}_{c},\, q\in
(0,1)$ with $Y$ of class (D) if $p=1$. We will show that
nevertheless there exists a solution with weaker integrability
properties. Before proving our main result let us note that in
\cite{EKPPQ,HP,C} reflected BSDEs with generator $f$ such that
$|f(t,y,z)|\le M(|f(t,0,0)|+|y|+|z|)$ for some $M\ge0$ are
considered. In case $p=2$ it is proved there that if we assume
that $\xi, \int_{0}^{T}|f(s,0,0)|\,ds\in \mathbb{L}^{2}(\FF_{T})$,
$L$ is continuous and $L^{+}\in \mathcal{S}^{2}$ then there exists
a solution $(Y,Z,K)\in \mathcal{S}^{2}\otimes M^{2}\otimes
\mathcal{V}^{+,2}_{c}$ of (\ref{eq1.1}) (see \cite{EKPPQ} for the
case of Lipschitz continuous generator and \cite{HP,C} for
continuous generator). We would like to stress that although in
\cite{EKPPQ,HP,C} condition (H7) is not explicitly stated, it is
satisfied, because if $f$ satisfies the linear growth condition
and $L^{+}\in\mathcal{S}^{2}$ then
\[
E(\int_{0}^{T}f^{-}(t,L^{+,*}_{t},0)\,dt)^{2} \le
2M^{2}T^{2}+2T^{2}E|L^{+,*}_{T}|^{2}<+\infty
\]
and $L_{t}\le L^{+,*}_{t}$, $t\in[0,T]$,
$L^{+,*}\in\mathcal{V}^{+,2}_{c}$.

\begin{tw}
\label{tw7.1} Let \mbox{\rm (H1)}--\mbox{\rm (H6)} (resp.
\mbox{\rm (H1)}--\mbox{\rm(H6)}, \mbox{\rm (Z)}) be satisfied and
$L^{+}\in\mathcal{S}^{p}$ for some $p>1$ (resp. $L^{+}$ is of
class \mbox{\rm(D)}). Then there exists a solution $(Y,Z,K)\in
\mathcal{S}^{p}\otimes M\otimes \mathcal{V}^{+}_{c}$ (resp.
$(Y,Z,K)\in \mathcal{S}^{q}\otimes M\otimes \mathcal{V}^{+}_{c}$,
$q\in (0,1)$ such that $Y$ is of class \mbox{\rm(D)}) of the
\mbox{\rm RBSDE}$(\xi,f,L)$.
\end{tw}
\begin{dow}
We first assume that $p=1$.  By \cite[Theorem 6.3]{BDHPS} there
exists a unique solution $(Y^{n},Z^{n})\in
\bigcap_{q<1}\mathcal{S}^{q}\otimes M^{q}$ of (\ref{eq5.01}) such
that $Y^n$  is of class (D). By Proposition \ref{tw2} (see also
the reasoning used at the beginning of the proof of Theorem
\ref{tw2}), for every $n\in\mathbb{N}$,  $Y^{n}_{t}\le
Y^{n+1}_{t}$ and $Y^{n}_{t}\le\bar{Y}^{n}_{t}$, $t\in [0,T]$,
where $(\bar{Y}^{n},\bar{Z}^{n})\in
\bigcap_{q<1}\mathcal{S}^{q}\otimes M^{q}$ is a solution of the
BSDE
\begin{align*}
\bar{Y}^{n}_{t}=\xi+\intt f^{+}(s,\bar{Y}^{n}_{s},
\bar{Z}^{n}_{s})\, ds+\intt n(\bar{Y}^{n}_{s}-L_{s})^{-}\,ds-\intt
\bar{Z}^{n}_{s}\,dB_{s},\quad t\in [0,T]
\end{align*}
such that $\bar{Y}^{n}$ is of class (D). Hence
\begin{align}
\label{apx1} |Y^{n}_{t}|\le |Y^{1}_{t}|+|\bar{Y}^{n}_{t}|,\quad
t\in [0,T].
\end{align}
Put
\[
R_{t}(L)=\esssup_{t\le\tau\le T} E(L_{\tau}|\FF_{t}).
\]
It is known (see \cite{CK,E}) that $R(L)$ has a continuous version
(still denoted by $R(L)$) such that $R(L)$ is a supermartingale of
class (D) majorizing the process $L$. Moreover, by the Doob-Meyer
decomposition theorem there exist a uniformly integrable
continuous martingale $M$ and a process
$V\in\mathcal{V}^{+,1}_{c}$ such that $R(L)=M+V$. In particular,
by \cite[Lemma 6.1]{BDHPS}, $R(L)\in
\MM^{q}_{c}+\mathcal{V}^{+,1}_{c}$ for every $q\in (0,1)$.
Therefore the data $\xi,f^{+},L$ satisfy assumptions (H1)--(H6),
(Z) and (H7*) with $X=R(L)$. Hence, by Theorem \ref{tw2}, there
exists a unique solution $(\bar{Y},\bar{Z},\bar{K})\in
\mathcal{S}^{q}\otimes M^{q}\otimes \mathcal{V}^{+,1}_{c}$, $q\in
(0,1)$, of the RBSDE$(\xi,f^{+},L)$ such that $\bar Y$ is of class
(D) and
\begin{align*}
\bar{Y}^{n}_{t}\nearrow \bar{Y}_{t},\quad t\in [0,T].
\end{align*}
By the above and (\ref{apx1}),
\begin{align}
\label{apx3} |Y^{n}_{t}|\le |Y^{1}_{t}|+|\bar{Y}_{t}|,\quad t\in
[0,T].
\end{align}
Put $Y_{t}=\sup_{n\ge1} Y^{n}_{t}$, $t\in [0,T]$ and
\[
\tau_{k}=\inf\{t\in [0,T]; \int_{0}^{t}f^{-}(s,R_{s}(L),0)\,ds
>k\}\wedge T.
\]
Then $f,L$ satisfy assumptions (H1)--(H6), (Z) and (H7*) with
$X=R(L)$ on each interval $[0,\tau_{k}]$. Therefore analysis
similar to that in the proof of (\ref{eq5.03}), but applied to the
equation
\begin{equation} \label{apx4}
Y^{n}_{t\wedge\tau_{k}}=Y^{n}_{\tau_{k}}
+\int_{t\wedge\tau_{k}}^{\tau_{k}}f(s,Y^{n}_{s},Z^{n}_{s})\,ds
+\int^{\tau_{k}}_{t\wedge
\tau_{k}}n(Y^{n}_{s}-L_{s})^{-}\,ds-\int_{t\wedge
\tau_{k}}^{\tau_{k}} Z^{n}_{s}\, dB_{s}
\end{equation}
instead of (\ref{eq5.01}), shows that for every $k\in\mathbb{N}$,
\begin{align}
\label{apx5} E\sup_{0\le t\le \tau_{k}}|Y^{n}_{t}-Y^{m}_{t}|^{q}
+E(\int_{0}^{\tau_{k}}|Z^{n}_{s}-Z^{m}_{s}|^{2}\,ds)^{q/2}
+E\sup_{0\le t\le\tau_{k}}|K^{n}_{t}-K^{m}_{t}|^{q}\rightarrow0
\end{align}
as $n,m\rightarrow +\infty$, where
$K^{n}_{t}=\int_{0}^{t}n(Y^{n}_{s}-L_{s})^{-}\,ds$. (The only
difference between the proof of (\ref{apx5}) and (\ref{eq5.03}) is
caused by the fact that in (\ref{apx4}) the terminal condition
$Y^{n}_{\tau_{k}}$ depends on $n$. But in view of (\ref{apx3}),
monotonicity of the sequence $\{Y^{n}\}$ and integrability of
$Y^{1},\bar{Y}$ the dependence of $Y^{n}_{\tau_{k}}$ on $n$
presents no difficulty). Since the sequence $\{\tau_{k}\}$ is
stationary, from (\ref{apx4}), (\ref{apx5}) we conclude that there
exist $K\in \mathcal{V}^{+}_{c}$ and $Z\in M$ such that
\[
Y_{t}=\xi+\intt f(s,Y_{s},Z_{s})\,ds+\intt dK_{s} -\intt
Z_{s}\,dB_{s},\quad t\in [0,T]
\]
and (\ref{apx5}) holds with $(Y,Z,K)$ in place of
$(Y^{n},Z^{n},K^{n})$. From the properties of the sequence
$\{(Y^{n},Z^{n},K^{n})\}$ on $[0,\tau_{k}]$ proved in Theorem
\ref{tw2} it follows that
\[
Y_{t}\ge L_{t},\quad t\in [0,\tau_{k}],
\quad \int_{0}^{\tau_{k}}(Y_{s}-L_{s})\,ds=0
\]
for $k\in \mathbb{N}.$ Due to stationarity of the sequence
$\{\tau_{k}\}$ this  implies that
\[
Y_{t}\ge L_{t},\quad t\in [0,T],\quad \int_{0}^{T}(Y_{s}-L_{s})\,ds=0.
\]
Accordingly, the triple $(Y,Z,K)$ is a solution of
RBSDE$(\xi,f,L)$.

In case $p>1$ the proof is similar. As a matter of fact it is
simpler because instead of considering the Snell envelope $R(L)$
of the process $L$ it suffices to consider the process $L^{+,*}$.
\end{dow}

\begin{uw}
From Proposition \ref{stw4} it follows that  the solution obtained
in Theorem \ref{tw7.1} is unique in its class for $p>1$.  In case
$p=1$ it is unique in its class if $f$ does not depend on $z$ (see
Remark \ref{uw.dic}).
\end{uw}

The next example shows that in general the process $K$ of Theorem
\ref{tw7.1} may be nonintegrable for any $q>0$.
\begin{prz}
Let $f(t,y)=-y^{+}\exp(|B_{t}|^{4})$, $L_{t}\equiv 1$, $\xi\equiv
1$. Then $\xi,f,L$ satisfy (H1)--(H6) and $L\in\mathcal{S}^{p}$
for every $p\ge 1$. So by Theorem \ref{tw7.1} and Proposition
\ref{stw2.5} there exists a unique solution $(Y,Z,K)\in
\mathcal{S}^{2}\otimes M\otimes \mathcal{V}^{+}_{c}$ of the
RBSDE$(\xi,f,L)$. Observe that $EK_{T}^{q}=+\infty$ for any $q>0$.
Of course, to check this it suffices to consider the case
$q\in(0,1]$. Aiming for a contradiction, suppose that $q\in (0,1]$
and $EK^{q}_{T}<+\infty$. Then  by \cite[Lemma 3.1]{BDHPS}, $Z\in
M^{q}$, which implies that $E(\int_{0}^{T}
f^{-}(t,Y_{t})\,dt)^{q}<+\infty$. On the other hand, since
$Y_{t}\ge 1$ for $t\in [0,T]$, it follows that
\[
E(\int_{0}^{T} f^{-}(t,Y_{t})\,dt)^{q}\ge
E\int_{0}^{T}(f^{-}(t,1))^{q}\,dt
=E\int_{0}^{T}\exp(q|B_{t}|^{4})\,dt=+\infty.
\]
\end{prz}

\end{document}